\documentclass[12pt]{article}

\usepackage{amsfonts}

\usepackage[utf8]{inputenc}

\usepackage[T1]{fontenc}

\usepackage{textalpha} 

\usepackage{bm}
\usepackage{enumerate} 
\usepackage{amssymb,amsmath}
\usepackage{caption}
\usepackage{subcaption}
\usepackage{accents} 

\usepackage{accents} 

\let\Horig\H

\usepackage{tikz}
\usetikzlibrary{automata,topaths}
\usetikzlibrary{shapes}
\usetikzlibrary{plotmarks}

\usepackage{float} 
\usepackage{textcomp} 
\usepackage{siunitx}
\usepackage{color} 
\definecolor{lightblue}{rgb}{0,0.2,0.5}
\usepackage[colorlinks=true, urlcolor=lightblue,linkcolor=lightblue, citecolor=lightblue]{hyperref}
\usepackage{ucs}

\definecolor{blue1}{rgb}{0,0.1,0.9}
\definecolor{mauve}{rgb}{0.7,0,0.43}
\definecolor{dkgreen}{rgb}{0,0.6,0}

\definecolor{codegreen}{rgb}{0,0.6,0}
\definecolor{codegray}{rgb}{0.5,0.5,0.5}
\definecolor{backcolour}{rgb}{0.97,0.97,0.97}

\usepackage{listingsutf8}
\lstdefinelanguage{Sage}[]{Python}
{morekeywords={False,sage,True},sensitive=true}

\definecolor{dblackcolor}{rgb}{0.0,0.0,0.0}
\definecolor{dbluecolor}{rgb}{0.01,0.02,0.7}
\definecolor{dgreencolor}{rgb}{0.2,0.4,0.0}
\definecolor{dgraycolor}{rgb}{0.30,0.3,0.30}

\usepackage{fvextra}

\usepackage{graphicx}
\usepackage{flushend,cuted}
\usepackage{bm}
\usepackage{tabularx}

\usepackage{indentfirst}
\usepackage{amssymb}
\usepackage{xparse}
\usepackage{tikz}
\usepackage{mdwlist}
\usepackage{tkz-graph}

\DeclareMathAlphabet{\eufrak}{U}{}{}{} 
\SetMathAlphabet\eufrak{normal}{U}{euf}{m}{n}
\SetMathAlphabet\eufrak{bold}{U}{euf}{b}{n}

\newtheorem{assumption}{Assumption}[section]

\newcommand{\R}{\mathbb{R}}

\newcommand{\E}{\mathbb{E}}
\newcommand{\IP}{\mathbb{P}}

\newcommand{\bone}{{\bf 1}}

\newcommand{\inte}{\mathbb{N}}

\makeatletter
\newcommand*\rel@kern[1]{\kern#1\dimexpr\macc@kerna}
\newcommand*\widebar[1]{%
  \begingroup
  \def\mathaccent##1##2{%
    \rel@kern{0.8}%
    \overline{\rel@kern{-0.8}\macc@nucleus\rel@kern{0.2}}%
    \rel@kern{-0.2}%
  }%
  \macc@depth\@ne
  \let\math@bgroup\@empty \let\math@egroup\macc@set@skewchar
  \mathsurround\z@ \frozen@everymath{\mathgroup\macc@group\relax}%
  \macc@set@skewchar\relax
  \let\mathaccentV\macc@nested@a
  \macc@nested@a\relax111{#1}%
  \endgroup
}
\makeatother

\makeatletter
\DeclareRobustCommand\widecheck[1]{{\mathpalette\@widecheck{#1}}}
\def\@widecheck#1#2{%
    \setbox\z@\hbox{\m@th$#1#2$}%
    \setbox\tw@\hbox{\m@th$#1%
       \widehat{%
          \vrule\@width\z@\@height\ht\z@
          \vrule\@height\z@\@width\wd\z@}$}%
    \dp\tw@-\ht\z@
    \@tempdima\ht\z@ \advance\@tempdima2\ht\tw@ \divide\@tempdima\thr@@
    \setbox\tw@\hbox{%
       \raise\@tempdima\hbox{\scalebox{1}[-1]{\lower\@tempdima\box
\tw@}}}%
    {\ooalign{\box\tw@ \cr \box\z@}}}
\makeatother

\newtheorem{prop}{Proposition}[section]
\newtheorem{lemma}[prop]{Lemma}
\newtheorem{definition}[prop]{Definition}
\newtheorem{corollary}[prop]{Corollary}
\newtheorem{thm}[prop]{Theorem}

\def\({\left(}
\def\){\right)}

\def\[{\left[}
\def\]{\right]}
\def\real{{\mathord{\mathbb R}}}

\def\P{\mathbb{P}}

\newcommand{\EE}{\mathsf{EP}}
\newcommand{\GG}{\mathsf{G}}

\newenvironment{Proof}{\removelastskip\par\medskip
\noindent{\em Proof.} \rm}{\penalty-20\null\hfill$\square$\par\medbreak}

\allowdisplaybreaks

\numberwithin{equation}{section}

\GraphInit[vstyle = Shade]
\usetikzlibrary[intersections,
positioning,
petri,
backgrounds,
fit,
decorations.pathmorphing,
arrows,
arrows.meta,
bending,
calc,
intersections,
through,
backgrounds,
shapes.geometric,
quotes,
matrix,
trees,
shapes.symbols,
graphs,
math,
patterns,
external,
scopes,
matrix,
lindenmayersystems,
shapes.callouts,
shapes.misc,
angles,
shapes.arrows,
shadings]
\tikzset{snake it/.style={-stealth,
decoration={snake, 
    amplitude = .4mm,
    segment length = 2mm,
    post length=0.9mm},decorate}}

\usetikzlibrary{matrix,calc}

\newcommand{\convexpath}[2]{
[   
    create hullnodes/.code={
        \global\edef\namelist{#1}
        \foreach [count=\counter] \nodename in \namelist {
            \global\edef\numberofnodes{\counter}
            \node at (\nodename) [draw=none,name=hullnode\counter] {};
        }
        \node at (hullnode\numberofnodes) [name=hullnode0,draw=none] {};
        \pgfmathtruncatemacro\lastnumber{\numberofnodes+1}
        \node at (hullnode1) [name=hullnode\lastnumber,draw=none] {};
    },
    create hullnodes
]
($(hullnode1)!#2!-90:(hullnode0)$)
\foreach [
    evaluate=\currentnode as \previousnode using \currentnode-1,
    evaluate=\currentnode as \nextnode using \currentnode+1
    ] \currentnode in {1,...,\numberofnodes} {
-- ($(hullnode\currentnode)!#2!-90:(hullnode\previousnode)$)
  let \p1 = ($(hullnode\currentnode)!#2!-90:(hullnode\previousnode) - (hullnode\currentnode)$),
    \n1 = {atan2(\y1,\x1)},
    \p2 = ($(hullnode\currentnode)!#2!90:(hullnode\nextnode) - (hullnode\currentnode)$),
    \n2 = {atan2(\y2,\x2)},
    \n{delta} = {-Mod(\n1-\n2,360)}
  in 
    {arc [start angle=\n1, delta angle=\n{delta}, radius=#2]}
}
-- cycle
}

\tikzset{hide labels/.style={every label/.append style={text opacity=0}}}

\makeatletter
\lst@InputCatcodes
\def\lst@DefEC{%
 \lst@CCECUse \lst@ProcessLetter
  ^^80^^81^^82^^83^^84^^85^^86^^87^^88^^89^^8a^^8b^^8c^^8d^^8e^^8f%
  ^^90^^91^^92^^93^^94^^95^^96^^97^^98^^99^^9a^^9b^^9c^^9d^^9e^^9f%
  ^^a0^^a1^^a2^^a3^^a4^^a5^^a6^^a7^^a8^^a9^^aa^^ab^^ac^^ad^^ae^^af%
  ^^b0^^b1^^b2^^b3^^b4^^b5^^b6^^b7^^b8^^b9^^ba^^bb^^bc^^bd^^be^^bf%
  ^^c0^^c1^^c2^^c3^^c4^^c5^^c6^^c7^^c8^^c9^^ca^^cb^^cc^^cd^^ce^^cf%
  ^^d0^^d1^^d2^^d3^^d4^^d5^^d6^^d7^^d8^^d9^^da^^db^^dc^^dd^^de^^df%
  ^^e0^^e1^^e2^^e3^^e4^^e5^^e6^^e7^^e8^^e9^^ea^^eb^^ec^^ed^^ee^^ef%
  ^^f0^^f1^^f2^^f3^^f4^^f5^^f6^^f7^^f8^^f9^^fa^^fb^^fc^^fd^^fe^^ff%
  ^^^^03b6
  ^^^^03b1^^^^03b2^^^^03b3%
  ^^00}
\lst@RestoreCatcodes
\makeatother

\begin{document}
\title{
\huge
 Graph connectivity with fixed endpoints in the random-connection model
} 

\author{
  Qingwei Liu\footnote{\href{mailto:xiaozong30@gmail.com}{xiaozong30@gmail.com}}
  \qquad
      Nicolas Privault\footnote{
\href{mailto:nprivault@ntu.edu.sg}{nprivault@ntu.edu.sg}
}
  \\
\small
Division of Mathematical Sciences
\\
\small
School of Physical and Mathematical Sciences
\\
\small
Nanyang Technological University
\\
\small
21 Nanyang Link, Singapore 637371
}

\maketitle

\vspace{-0.5cm}

\begin{abstract} 
We consider the count of subgraphs with an arbitrary configuration of endpoints in the random-connection model based on a Poisson point process on $\real^d$. We present combinatorial expressions for the computation of the cumulants and moments of all orders of such subgraph counts, which allow us to estimate the growth of cumulants as the intensity of the underlying Poisson point process goes to infinity. As a consequence, we obtain a central limit theorem with explicit convergence rates under the Kolmogorov distance, and connectivity bounds. Numerical examples are presented using a computer code in SageMath for the closed-form computation of cumulants of any order, for any type of connected subgraph, and for any configuration of endpoints in any dimension $d\geq 1$. In particular, graph connectivity estimates, Gram-Charlier expansions for density estimation, and correlation estimates for joint subgraph counting are obtained. 
\end{abstract}


\baselineskip0.7cm

\section{Introduction}
\noindent
This paper considers the statistics and asymptotic behavior of
subgraph counts in a multidimensional random-connection
model based on a Poisson point process, 
 which can be used to model physical systems in e.g.
 statistical mechanics 
 \cite{kartungiles2016}, 
 wireless networks
 \cite{ta2007,Mao2010,georgiou2015}, 
 or cosmology \cite{cunningham2017,fountoulakis2020}. 

 \medskip
 
 The random-connection model,
 in which vertices are randomly located and connected
 with location-dependent probabilities, 
 is a natural generalization of
 e.g. the Erd\H os-R\'enyi random graph
 or the stochastic block model \cite{snijders}.  
 Namely, given $\mu$ a 
 diffuse Radon measure on $\R^d$, the random-connection model $G_H (\eta )$
 consists of an underlying Poisson point process $\eta$ on $\R^d$
 with intensity of the form $\lambda \mu(\mathrm{d}x)$, $\lambda >0$,  
 in which any two vertices $x,y$ in $\eta $ are connected
 with the probability $H(x,y)$, 
 where $H:\real^d \times \real^d \to [0,1]$ is a
 {symmetric} connection function. 
 
\medskip 

{
  In addition to modeling the
  random locations of network nodes,
  many applications of wireless networks
  require the use of endpoints which are physical devices
  placed at given fixed locations,
  such as for example roadside units
  in vehicular networks such as VANETs,
  see, e.g., 
  \cite{ng2011}, \cite{zhang2012}.}
The count of subgraphs that connect any single point 
$x$ in the Poisson process $\eta$
to $m$ fixed endpoints $y_1,\ldots , y_m\in \real^d$
 is known to have a Poisson distribution with mean
$
 \lambda \int_{\real^d} H(x,y_1)\cdots H(x,y_m) \ \! \mu (\mathrm{d}x)$,
 see e.g. \S~4 in \cite{prkhp}.
 This Poisson property has been used in \cite{giles-privault2_published} 
 to derive closed-form estimates of two-hop
 connectivity in the random-connection model
 when $m=2$, see Proposition~III.2 therein.

\medskip 

In this paper, we consider the count of general connected
 subgraphs with a general configuration of
 fixed endpoints at fixed locations $y_1 , \ldots,y_m \in \R^d$ 
 in the random-connection model  
 $G_H (\eta \cup \{y_1, \ldots ,y_m\})$ constructed on the union
 of the Poisson point process $\eta$ and $\{ y_1, \ldots ,y_m\}$.
 {
   In particular, we extend the subgraph count cumulant formulas obtained
   on $G_H (\eta )$ in \cite{LiuPrivault}  
   by taking into account the presence of endpoints in
   $G_H (\eta \cup \{y_1, \ldots ,y_m\})$,
   and we provide
   SageMath coding implementations
   for joint cumulant expressions of any order.} 
 \label{here}
 
 \medskip 

In Proposition~\ref{mom-cumfor} we derive general expressions
for the moments and cumulants of the count $N^G_{y_1,\ldots , y_m}$ of subgraphs with fixed endpoints
$y_1, \ldots ,y_m$ in $G_H (\eta\cup \{y_1, \ldots ,y_m\})$.
Such expressions allow us to determine the dominant terms in the growth of
 cumulants as the intensity $\lambda$ of the underlying point process tends to infinity, 
 by estimating the counts of vertices and edges in connected
 partition diagrams as in e.g. \cite{khorunzhiy}. 
 As a consequence, in Theorem~\ref{khopone}
 we obtain growth estimates for the cumulants of
 the subgraph count $N^G_{y_1,\ldots , y_m}$. 

 \medskip

 This allows us to show the convergence of renormalized subgraph
 counts to the normal distribution in Proposition~\ref{fjlfa12}
 as the intensity $\lambda$ of the underlying Poisson point process on $\R^d$
 tends to infinity.
Convergence rates under the Kolmogorov distance
are then obtained in Proposition~\ref{pkol}
for the normal approximation of subgraph counts
from the combinatorics of cumulants
 and the {Statulevi\v{c}ius condition}, see \cite{rudzkis,doering}
 and Lemma~\ref{Statuleviciuscond1}, 
 extending the results obtained in \cite{LiuPrivault}
 for subgraphs without endpoints.
 See also \cite{thale18} 
 for other applications of this condition
 to concentration inequalities,
 normal approximation and moderate deviations for random polytopes.
 In Proposition~\ref{jklf3}, connectivity probability
 estimates and bounds are derived using the second moment method 
 and the factorial moment expansions in Proposition~\ref{fdshkf0}. 
 
\medskip 

In Section~\ref{examples} we consider several examples of
subgraphs with endpoints such as $k$-hop paths,
triangles and trees, for which  
exact cumulant computations are matched to their
Monte Carlo estimates 
 using the Rayleigh connection function $H(x,y) = e^{ - \beta \Vert x - y\Vert^2}$, 
 $\beta > 0$.
 In those examples
 we obtain graph connectivity estimates,
 and correlation estimates for joint graph counting,
 which are matched to the outputs of Monte Carlo simulations. 
 In addition, using third order cumulant expressions, we also provide
 improved fits of probability density functions of renormalized
 subgraph counts when the Gaussian approximation is not valid,
 see Figure~\ref{fig5}. 

 \medskip
 
 Computations are done in closed form using symbolic calculus 
 in the SageMath coding implementations 
 presented in Appendices~\ref{fjkldsf}-\ref{fjkldsf-2}, 
 and available for download at 
 \url{https://github.com/nprivaul/random-connection}.  
 We note that although intensive computations may be required,
the types of connected subgraphs and associated configurations 
of endpoints considered is only limited by the available computing
power. 

\medskip 

This paper is organised as follows.
 Section~\ref{rcm} introduces some 
 preliminaries on subgraph counting
 and the computation of moments
 using summations over partitions
 in the random-connection model. 
 In Section~\ref{diagramrepresentation}, we use
 partition diagrams to compute the cumulants
 of the counts of subgraphs with endpoints in the random-connection
model. 
Subgraph count asymptotics and the associated
central limit theorem are given 
in Section~\ref{sca}, and numerical examples are
presented in Section~\ref{examples}.
A general derivation of joint cumulant identities is given 
in Appendix~\ref{appendixa}, extending the construction of
\cite{LiuPrivault} from the univariate to the multivariate case,
for use in Section~\ref{corrof}. 
 Basic results on Gram-Charlier expansions and probability approximation
 using cumulant and moment methods are recalled in 
 Appendices~\ref{statuleviciuscond} and \ref{s5}.
 The SageMath codes for the computation of cumulants and joint cumulants are listed in Appendices~\ref{fjkldsf} and \ref{fjkldsf-2},
 and available for download at 
\url{https://github.com/nprivaul/random-connection}.  

\section{Subgraph counts in the random-connection model} 
\label{rcm}
\begin{figure}[H]
  \centering
  \begin{tikzpicture}
\draw[black, thick] (-1,1) rectangle (6,3);
\filldraw [gray] (1,2) circle (2pt);
\node[font=\small] at (1,2.4) {$v_1$};
\filldraw [gray] (2,2) circle (2pt);
\node[font=\small] at (2,2.4) {$v_2$};
\filldraw [gray] (3,2) circle (2pt);
\node[font=\small] at (3,2.4) {$v_3$};
\filldraw [gray] (4,2) circle (2pt);
\node[font=\small] at (4,2.4) {$v_4$};
\filldraw [gray] (-0,2) circle (2pt);
\node[font=\small] at (-0.5,2) {$w_1$};
\filldraw [gray] (5,2) circle (2pt);
\node[font=\small] at (5.5,2) {$w_2$};
\draw[thick,blue] (1,2) .. controls (1.5,2.3) .. (2,2);
\draw[thick,blue] (2,2) .. controls (2.5,2.3) .. (3,2);
\draw[thick,blue] (2,2) .. controls (3,1.7) .. (4,2);
\draw[thick,blue] (1,2) -- (0,2);
\draw[thick,blue] (4,2) -- (5,2);
 \begin{pgfonlayer}{background}
    \filldraw [line width=4mm,black!3]
      (-0.8,1.2)  rectangle (5.8,2.8);
  \end{pgfonlayer}
\end{tikzpicture}
\caption{
 Graph $G=(V_G,E_G)$ with $V_G=(v_1, v_2,v_3,v_4;w_1,w_2)$, 
 $n=3$, $r=4$, $m=2$.}
\label{fig:diagram0}
\end{figure}
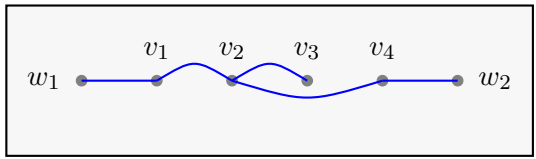

\noindent 
 In what follows we consider 
 a Radon measure $\mu$ on $\real^d$, and we let 
 $\IP_\lambda$, $\lambda > 0$, 
 denote the distribution of the Poisson point process $\eta$ 
 with intensity $\lambda \mu (\mathrm{d}x)$ on the space 
 $$\mathcal{C}:=\left\{\eta \subset \R^d \ : \ |\eta \cap A|<\infty ~\text{for any bounded set $A\subset \R^d$}\right\}  
 $$
 of locally finite configurations on $\R^d$, whose elements 
 $\eta \in \Omega$ are identified with the Radon point measures 
 $\displaystyle \eta = \sum_{x\in \eta } \epsilon_x$, 
 so that $\eta (B)$ represents the random
 number of points contained in a Borel set in $\R^d$. 
 In other words, 
 \begin{enumerate}[i)] 
 \item for any relatively compact Borel set $B\subset \R^d$, the distribution of $\eta(B)$ under $\IP_\lambda$ is Poisson with parameter $\lambda \mu (B)$; 

\vspace{-0.2cm}  

\item for any $n\geq 2$ and pairwise disjoint
  relatively compact Borel sets $B_1, \ldots ,B_n\subset \R^d$, the random
  variables $\eta(B_1), \ldots ,\eta(B_n)$ are independent under $\IP_\lambda$.
\end{enumerate}
\noindent
 For $n\geq 1$ we let $[n]:=\{1, \ldots ,n\}$,
{where $n$ will later on denote the order of
 the considered moments and cumulants of subgraph counts,} 
 and for any set $A$ we denote by $\Pi (A)$
 the collection of all set partitions of $A$.
 We also let $|A|$ denote the number of elements of any finite set $A$,
 and, in particular, $|\sigma|$ represents the number of blocks in a partition
 $\sigma\in\Pi ([n]\times[r])$.
{Our approach to the computation of moments relies on 
 moment identities on the following form,
 see Proposition~3.1 in \cite{momentpoi}, and 
 Proposition~\ref{moment-1} for its multivariate generalization.
\begin{prop}
\label{moment-1-0}
 Let $n\ge1$ and $r \ge1$,
 and let $f :(\R^d)^r \to\R$ be a sufficiently integrable
 measurable function. We have
$$
  \E\left[\left( 
 \sum_{(x_1, \ldots ,x_r)\in\eta^r } f (x_1,\dots,x_r) 
 \right)^n
  \right]
  =
  \hskip-0.3cm
  \sum_{\rho\in\Pi([n]\times[r])}
  \hskip-0.3cm
  \lambda^{|\rho|}
  \hskip-0.1cm
  \int_{(\R^d)^{|\rho|}}
 \prod_{k=1}^n
 f\big(
 x_{\zeta_\rho (k,1)}, \ldots , x_{\zeta_\rho (k,r)}
 \big) \ \! 
 \mu(\mathrm{d}x_1) \cdots \mu(\mathrm{d}x_{|\rho|}), 
$$
 where, for $\rho = \{ \rho_1,\ldots , \rho_{|\rho |} \}$ a partition of
 $[n]\times[r]$, we let
 $\zeta_\rho (k,l)$ denote the index $p$ of the block 
 $\rho_p$ of $\rho$ to which $(k,l)$ belongs. 
\end{prop}
 In particular, Proposition~\ref{moment-1-0}
 will yield cumulant expressions from M\"obius inversion and 
 combinatorial arguments based on 
\cite{MalyshevMinlos91}, 
\cite{khorunzhiy} and 
\cite{LiuPrivault},
see Propositions~\ref{mom-cumfor} and \ref{fjklf2}.}
\begin{definition} 
 Given $H:\R^d\times \R^d\to[0,1]$ a symmetric connection function and  
 $y_1 , \ldots,y_m$ fixed points in $\R^d$,
 the random-connection model 
  $G_H (\eta \cup \{y_1, \ldots ,y_m\})$
 is the random graph built on the union of $\{y_1 , \ldots,y_m \}$
 and a Poisson point process
 sample $\eta$, 
  in which any two distinct points 
  $x,y\in \eta \cup \{y_1, \ldots ,y_m\}$
  are independently connected by an edge with the probability $H(x,y)$.  
\end{definition} 
{In the sequel we will consider
  a family of connected graphs with endpoints
  which are described in the following assumption.
}
\begin{assumption}
   \label{a0}
 Given $r\geq 2$ and $m\geq 0$, we 
consider a connected graph $G=(V_G,E_G)$ with edge set $E_G$ and
vertex set 
$V_G=(v_1, \ldots ,v_r; {w_1,\ldots , w_m})$, such that
\begin{enumerate}[i)]
\item the subgraph $\GG$ induced by $G$ on $\{v_1, \ldots ,v_r\}$ is connected, and 
\item \label{ii}
  {the endpoint vertices $w_1, \ldots ,w_m$} are not adjacent to each other in $G$. 
\end{enumerate}
In case $m=0$,
 Condition~\eqref{ii} is void and
$V_G=(v_1, \ldots ,v_r)$.  
\end{assumption}
In Figure~\ref{fig:diagram0} 
 an example of a graph {satisfying Assumption~\ref{a0}} 
 is described with $r=4$ and $m=2$. 

\vspace{-0.3cm}
\noindent
{As a convention, in the next definition
 the sets $\{w_1,\ldots , w_m\}$ and 
 $\{ y_1 , \ldots,y_m \} \subset \R^d$ are empty when $m=0$.} 
\begin{definition}
\label{fjkl} 
{Let $G$ be
a graph satisfying Assumption~\ref{a0}.} 
Given $m \geq 0$ fixed points $y_1 , \ldots,y_m \in \R^d$,
  for a.s. $\eta$ we let $N_{y_1,\ldots , y_m}^G$ denote the count of subgraphs
 in $G_H (\eta \cup \{y_1,\ldots , y_m \} )$
 that are isomorphic to $G=(V_G,E_G)$ in the sense that 
 there exists a (random) injection
 from $V_G$ into $\eta \cup \{y_1,\ldots , y_m \}$
 which is one-to-one from $\{{w_1,\ldots , w_m}\}$ to $\{y_1,\ldots , y_m\}$, 
 and preserves the graph structure of $G$. 
\end{definition}
\noindent
According to Definition~\ref{fjkl}, we express the subgraph count
 $N^G_{y_1,\ldots , y_m}$ as  
\begin{equation}
\nonumber
  N_{y_1,\ldots , y_m}^G=\sum_{(x_1, \ldots ,x_r)\in\eta^r }f_{y_1,\ldots , y_m} (x_1, \ldots ,x_r), 
\end{equation}
where the random
function $f:(\real^d)^r \to \{0,1\}$ defined as 
\begin{equation}
\nonumber
f_{y_1,\ldots , y_m} (x_1, \ldots ,x_r):=
\prod_{
  \substack{
    1 \leq i \leq r
    \\
    1 \leq j \leq m
    \\ \{v_i,w_j\}\in E_G }
}
\bone_{\{y_j\leftrightarrow x_i\}} 
\prod_{\substack{ 1 \leq k,l \leq r
    \\ \{v_k,v_l\}\in E_G}}\bone_{\{x_k \leftrightarrow x_l \}},
\qquad
 x_1,\ldots , x_r \in \R^d, 
\end{equation}
is independent of the Poisson point process $\eta$,
 and $\bone_{\{x\leftrightarrow y\}}=1$ if and only if
$x\neq y$ and $x,y\in \real^d$ are connected in the
 random-connection model 
 $G_H (\eta \cup \{y_1, \ldots ,y_m\})$.

\section{Partition diagrams} 
\label{diagramrepresentation}
\noindent 
{This section introduces the combinatorial background
needed for the derivation of moment and cumulant expressions
of subgraph counts.} 
 {The next definition introduces
   a notion of connectedness over the rows of 
   partitions of $[n]\times [r]$,
   and a flatness property which is satisfied when
   two indices on a same row belong to a given block,
   see Chapter~4 of \cite{peccatitaqqu}
   and Figure~\ref{fig:diagram2-0} below.} 
\begin{definition}
   \label{def-1}
   Given $n,r\geq 1$,
  let $\pi:=\{\pi_1, \ldots ,\pi_n\}$ be the partition in $\Pi ([n]\times[r])$
  given by 
  $$
  \pi_i:=\left\{(i,1), \ldots ,(i,r)\right\},
  \quad
  i=1, \ldots , n.
  $$
 \begin{enumerate}[i)]
   \item A set partition $\sigma\in\Pi ([n]\times[r])$ is connected if $\sigma\vee\pi=\widehat{1}$, 
     where
     $\sigma \vee\pi$ is the finest set partition which is coarser than both
     $\sigma$ and $\pi$, and $\widehat{1} = \{ [n]\times [r] \}$
is the coarsest partition of $[n]\times [r]$. 
\item 
 A set partition $\sigma\in\Pi ([n]\times[r])$ is non-flat if $\sigma\wedge\pi=\widehat{0}$,
 where
 $\sigma \wedge\pi$ is
 the coarsest set partition which is finer than both $\sigma$ and $\pi$,
 and $\widehat{0}$ is the finest partition of $[n]\times [r]$.
\end{enumerate} 
 We let $\Pi_{\widehat{1}} ([n]\times[r])$ denote the collection of all
connected partitions of $[n] \times [r]$. 
\end{definition}

\noindent 
 In the sequel, every partition
 $\rho \in \Pi(\pi_1\cup \cdots \cup \pi_n )$
 will be arranged into a 
 diagram denoted by $\Gamma(\rho ,\pi)$, 
 by arranging $\pi_1,\dots,\pi_n$ into $n$ rows 
 and connecting together the elements of every block of $\rho$. 
 Figure~\ref{fig:diagram2-0} presents two illustrations
 of flat non-connected and connected non-flat partition diagrams
 with $n=5$ and $r=4$, 
 {in which the partition $\rho$ is represented using
   line segments.} 
  
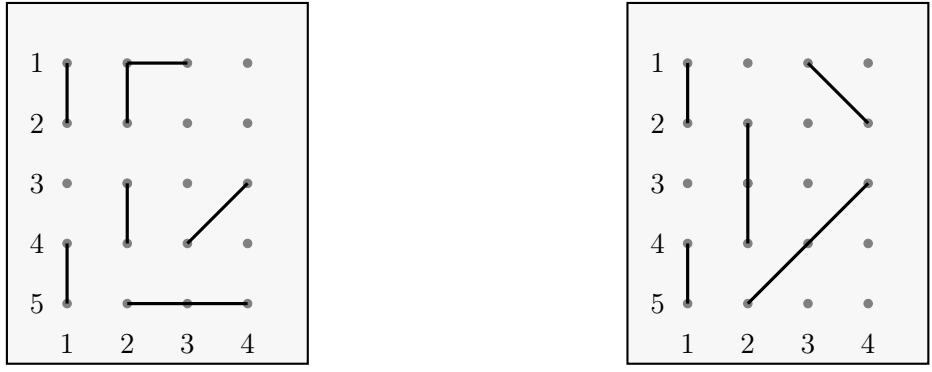
\begin{figure}[H]
\captionsetup[subfigure]{font=footnotesize}
\centering
\subcaptionbox{Flat non-connected diagram $\Gamma(\rho,\pi)$.}[.5\textwidth]{%
\begin{tikzpicture}[scale=0.8] 
\draw[black, thick] (0,0) rectangle (5,6);

\node[anchor=east,font=\small] at (0.8,5) {1};
\node[anchor=east,font=\small] at (0.8,4) {2};
\node[anchor=east,font=\small] at (0.8,3) {3};
\node[anchor=east,font=\small] at (0.8,2) {4};
\node[anchor=east,font=\small] at (0.8,1) {5};

\node[anchor=south,font=\small] at (1,0) {1};
\node[anchor=south,font=\small] at (2,0) {2};
\node[anchor=south,font=\small] at (3,0) {3};
\node[anchor=south,font=\small] at (4,0) {4};

\filldraw [gray] (1,1) circle (2pt);
\filldraw [gray] (2,1) circle (2pt);
\filldraw [gray] (3,1) circle (2pt);
\filldraw [gray] (4,1) circle (2pt);
\filldraw [gray] (1,2) circle (2pt);
\filldraw [gray] (2,2) circle (2pt);
\filldraw [gray] (3,2) circle (2pt);
\filldraw [gray] (4,2) circle (2pt);
\filldraw [gray] (1,3) circle (2pt);
\filldraw [gray] (2,3) circle (2pt);
\filldraw [gray] (3,3) circle (2pt);
\filldraw [gray] (4,3) circle (2pt);
\filldraw [gray] (2,3) circle (2pt);
\filldraw [gray] (1,4) circle (2pt);
\filldraw [gray] (2,4) circle (2pt);
\filldraw [gray] (3,4) circle (2pt);
\filldraw [gray] (4,4) circle (2pt);
\filldraw [gray] (1,5) circle (2pt);
\filldraw [gray] (2,5) circle (2pt);
\filldraw [gray] (3,5) circle (2pt);
\filldraw [gray] (4,5) circle (2pt);

\draw[very thick] (1,5) -- (1,4);
\draw[very thick] (2,5) -- (2,4);
\draw[very thick] (2,5) -- (3,5);

\draw[very thick] (1,2) -- (1,1);
\draw[very thick] (2,3) -- (2,2);
\draw[very thick] (2,1) -- (3,1) -- (4,1);
\draw[very thick] (3,2) -- (4,3);

 \begin{pgfonlayer}{background}
    \filldraw [line width=4mm,black!3]
      (0.27,0.27)  rectangle (4.76,5.76);
  \end{pgfonlayer}
\end{tikzpicture}}%
\subcaptionbox{Connected non-flat diagram $\Gamma(\rho,\pi)$.}[.5\textwidth]{
\begin{tikzpicture}[scale=0.8] 
\draw[black, thick] (0,0) rectangle (5,6);

\node[anchor=east,font=\small] at (0.8,5) {1};
\node[anchor=east,font=\small] at (0.8,4) {2};
\node[anchor=east,font=\small] at (0.8,3) {3};
\node[anchor=east,font=\small] at (0.8,2) {4};
\node[anchor=east,font=\small] at (0.8,1) {5};

\node[anchor=south,font=\small] at (1,0) {1};
\node[anchor=south,font=\small] at (2,0) {2};
\node[anchor=south,font=\small] at (3,0) {3};
\node[anchor=south,font=\small] at (4,0) {4};

\filldraw [gray] (1,1) circle (2pt);
\filldraw [gray] (2,1) circle (2pt);
\filldraw [gray] (3,1) circle (2pt);
\filldraw [gray] (4,1) circle (2pt);
\filldraw [gray] (1,2) circle (2pt);
\filldraw [gray] (2,2) circle (2pt);
\filldraw [gray] (3,2) circle (2pt);
\filldraw [gray] (4,2) circle (2pt);
\filldraw [gray] (1,3) circle (2pt);
\filldraw [gray] (2,3) circle (2pt);
\filldraw [gray] (3,3) circle (2pt);
\filldraw [gray] (4,3) circle (2pt);
\filldraw [gray] (2,3) circle (2pt);
\filldraw [gray] (1,4) circle (2pt);
\filldraw [gray] (2,4) circle (2pt);
\filldraw [gray] (3,4) circle (2pt);
\filldraw [gray] (4,4) circle (2pt);
\filldraw [gray] (1,5) circle (2pt);
\filldraw [gray] (2,5) circle (2pt);
\filldraw [gray] (3,5) circle (2pt);
\filldraw [gray] (4,5) circle (2pt);

\draw[very thick] (1,5) -- (1,4); 
\draw[very thick] (3,5) -- (4,4);

\draw[very thick] (1,2) -- (1,1);
\draw[very thick] (2,2) -- (2,4);
\draw[very thick] (2,1) -- (3,2) -- (4,3);

 \begin{pgfonlayer}{background}
    \filldraw [line width=4mm,black!3]
      (0.27,0.27)  rectangle (4.76,5.76);
  \end{pgfonlayer}
\end{tikzpicture}}%
\caption{Two examples of partition diagrams with $n=5$ and $r=4$.}
\label{fig:diagram2-0}
\end{figure}
\vspace{-0.4cm}

\begin{figure}[H]
\captionsetup[subfigure]{font=footnotesize}
\centering
\subcaptionbox{Diagram before merging edges and vertices.}[.5\textwidth]{%
\begin{tikzpicture}
\draw[black, thick] (0,0) rectangle (7,4);
\foreach \i in {1,2,3}
{
\filldraw [gray] (2,\i) circle (2pt);
\filldraw [gray] (3,\i) circle (2pt);
\filldraw [gray] (4,\i) circle (2pt);
\filldraw [gray] (5,\i) circle (2pt);
\draw[thick, dash dot,blue] (2,\i) .. controls (2.5,\i) .. (3,\i);
\draw[thick, dash dot,blue] (4,\i) .. controls (3.5,\i) .. (3,\i);
\draw[thick, dash dot,blue] (1,2) .. controls (1.5,\i) .. (2,\i);
\draw[thick, dash dot,blue] (6,2) .. controls (5.5,\i) .. (5,\i);
}

\draw[thick, dash dot,blue] (3,1) .. controls (4,1+.4) .. (5,1);
\draw[thick, dash dot,blue] (3,2) .. controls (4,2-.4) .. (5,2);
\draw[thick, dash dot,blue] (3,3) .. controls (4,3-.4) .. (5,3);

\node[anchor=north,font=\tiny] at (2,1) {(3,1)};
\node[anchor=north,font=\tiny] at (3,1) {(3,2)};
\node[anchor=north,font=\tiny] at (4,1) {(3,3)};
\node[anchor=north,font=\tiny] at (5,1) {(3,4)};
\node[anchor=south,font=\tiny] at (2,3) {(1,1)};
\node[anchor=south,font=\tiny] at (3,3) {(1,2)};
\node[anchor=south,font=\tiny] at (4,3) {(1,3)};
\node[anchor=south,font=\tiny] at (5,3) {(1,4)};
\node[anchor=south,font=\tiny] at (2,2) {(2,1)};
\node[anchor=south,font=\tiny] at (3,2) {(2,2)};
\node[anchor=south,font=\tiny] at (4,2) {(2,3)};
\node[anchor=south,font=\tiny] at (5,2) {(2,4)};

\filldraw [gray] (1,2) circle (2pt);
\node[anchor=east,font=\tiny] at (1,2) {$(1)$};
\filldraw [gray] (6,2) circle (2pt);
\node[anchor=west,font=\tiny] at (6,2) {$(2)$};
\draw[thick] (3,3) -- (3,2);
\draw[thick] (2,2) -- (2,1);
\draw[thick] (5,2) -- (5,1);
 \begin{pgfonlayer}{background}
    \filldraw [line width=4mm,black!3]
      (0.2,0.2)  rectangle (6.8,3.8);
  \end{pgfonlayer}
\end{tikzpicture}}%
\subcaptionbox{Graph $\rho_G$ after merging edges and vertices.}[.5\textwidth]{
\begin{tikzpicture}
\draw[black, thick] (0,0) rectangle (7,4);
\foreach \i in {3}
{
\filldraw [gray] (2,\i) circle (2pt);
\filldraw [gray] (3,\i) circle (2pt);
\filldraw [gray] (4,\i) circle (2pt);
\filldraw [gray] (5,\i) circle (2pt);
\draw[thick,blue] (2,\i) .. controls (2.5,\i) .. (3,\i);
\draw[thick,blue] (4,\i) .. controls (3.5,\i) .. (3,\i);
\draw[thick,blue] (3,\i) .. controls (4,\i-.4) .. (5,\i);
\draw[thick,blue] (1,2) .. controls (1.5,\i) .. (2,\i);
\draw[thick,blue] (6,2) .. controls (5.5,\i) .. (5,\i);
}
\node[anchor=south,font=\tiny] at (2,3) {1};
\node[anchor=south,font=\tiny] at (3,3) {2};
\node[anchor=south,font=\tiny] at (4,3) {3};
\node[anchor=south,font=\tiny] at (5,3) {4};
\node[anchor=south,font=\tiny] at (2,2) {5};
\node[anchor=south,font=\tiny] at (4,2) {6};
\node[anchor=south,font=\tiny] at (5,2) {7};
\node[anchor=north,font=\tiny] at (3,1) {8};
\node[anchor=north,font=\tiny] at (4,1) {9};
\node[anchor=north,font=\tiny] at (1,2) {10};
\node[anchor=north,font=\tiny] at (6,2) {11};
\filldraw [gray] (1,2) circle (2pt);
\filldraw [gray] (6,2) circle (2pt);
\filldraw [gray] (2,2) circle (2pt);
\filldraw [gray] (4,2) circle (2pt);
\filldraw [gray] (5,2) circle (2pt);
\filldraw [gray] (3,1) circle (2pt);
\filldraw [gray] (4,1) circle (2pt);
\draw[thick,blue] (1,2) .. controls (1.5,2) .. (2,2);
\draw[thick,blue] (3,3) .. controls (2.5,2.5) .. (2,2);
\draw[thick,blue] (3,3) .. controls (3.5,2.5) .. (4,2);
\draw[thick,blue] (3,3) .. controls (4,2.5) .. (5,2);
\draw[thick,blue] (5,2) -- (6,2);
\draw[thick,blue] (2,2) -- (3,1) -- (4,1);
\draw[thick,blue] (3,1) -- (5,2);

 \begin{pgfonlayer}{background}
    \filldraw [line width=4mm,black!3]
      (0.2,0.2)  rectangle (6.8,3.8);
  \end{pgfonlayer}
\end{tikzpicture}}%
\caption{
  Example of graph $\rho_G$
   with $n=3$, $r=4$, and $m=2$.}
\label{fig:diagram1}
\end{figure}
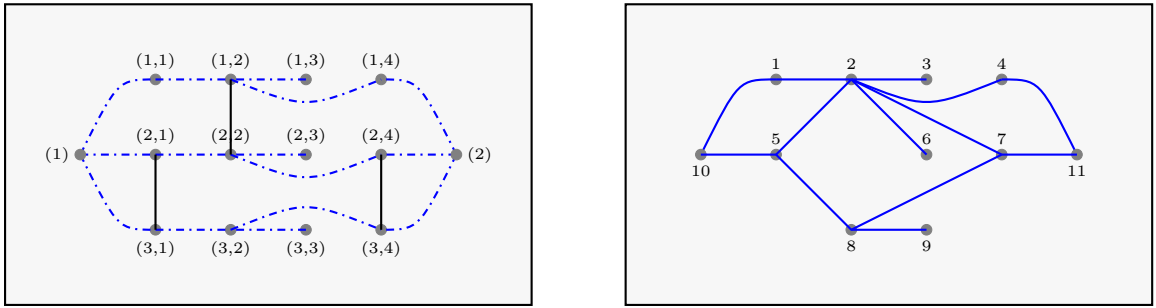

\noindent 
In Definition~\ref{fjklf},
 to any graph $G$ and set partition $\rho\in\Pi ([n]\times[r])$
 we associate a graph $\rho_G$ whose vertices are the blocks of $\rho$. {For this, we use 
   $n$ copies of the graphs induced by the $v_i$'s
   with addition of the end-points $w_1,\ldots , w_m$,
   and we merge the nodes obtained in this way
   on $[n]\times [r]$ according to the partition $\rho$.} 
\begin{definition}
   \label{fjklf}
   Given $\rho$ a partition of $[n]\times[r]$ 
   and $G=(V_G,E_G)$ a connected graph 
   on $V_G=(v_1, \ldots ,v_r; {w_1,\ldots , w_m})$, 
   we let $\rho_G$ denote the graph 
     constructed as follows on $[m] \cup [n]\times [r]$:
\begin{enumerate}[i)]  
\item for all $j_1, j_2\in [r]$, $j_1\not= j_2$, and $i\in [n]$, 
  an edge links $(i,j_1)$ to $(i,j_2)$ iff $\{v_{j_1},v_{j_2}\}\in E_G$; 
\item for all $(j,k)\in [r]\times [m]$ and $i\in [n]$, an edge
  links $(k)$ to $(i,j)$ iff $\{v_j,w_k\}\in E_G$; 
\item for all $i_1,i_2\in [n]$
  and $j_1,j_2\in [r]$,
  merge any two nodes $(i_1,j_1)$ and $(i_2,j_2)$ 
  if they belong to a same block in $\rho$;   
\item eliminating any redundant edges created by the above construction.
\end{enumerate}
\end{definition}
\noindent
 If $\rho\in\Pi ([n]\times[r])$
 takes the form $\rho = \{ b_1,\ldots , b_{|\rho |}\}$, 
 the graph $\rho_G$ forms a connected graph with
 $|\rho | + m$ vertices, and we reindex the set of vertices $V_{\rho_G}$
 of $\rho_G$  
 as $V_{\rho_G}=[|\rho | + m ]$ according to the lexicographic order
 on $\inte \times \inte$, 
 followed by the remaining $m$ vertices, 
 indexed as $\{|\rho |+1,\ldots , |\rho | +m\}$,
 see Figure~\ref{fig:diagram1}-$b)$
 in which we have $|\rho | =9$, $m = 2$, and $V_{\rho_G}=(1,\dots ,9;10, 11)$. 

 \medskip

 \noindent
{\bf Example}. 
Take $r=4$, $m=2$ and $V_G=(v_1, v_2,v_3,v_4;w_1,w_2)$. 
 Figure~\ref{fig:diagram1}-$b)$  
 shows the graph $\rho_G$ defined from 
 $G=(V_G,E_G)$ of Figure~\ref{fig:diagram0}
 and the $9$-block partition $\rho \in \Pi ([3]\times[4])$
 given by 
 \begin{align*}
   \rho = \big\{ & \{(1,1)\},
   \\
   & \{(1,2),(2,2)\},
      \\
   & \{(1,3)\},
   \\
   & \{(1,4)\},
   \\
   & \{(2,1),(3,1)\},
   \\
   & \{(2,3)\},
   \\
   & \{(2,4),(3,4)\},
   \\
   & \{(3,2)\},
   \\
   & \{(3,3)\}\big\}. 
\end{align*} 
 
  \noindent
 {
   In Figure~\ref{fig:diagram1}-$(a)$ 
   the partition $\rho$ is represented using line segments.
   } 

\vspace{-0.4cm}

\begin{definition}
  For $\rho\in\Pi ([n]\times[r])$ of the form
  $\rho = \{ b_1,\ldots , b_{|\rho |}\}$ 
  and $j \in [m]$, we let 
\begin{equation}
\nonumber
    {\cal A}^\rho_j:=\{ k \in [ |\rho | ] \ : \ \exists (s,i)\in b_k ~\mathrm{s.t.}~
    (v_i,w_j) \in E_G 
    \} 
\end{equation} 
denote the neighborhood of the vertex $(|\rho | + j)$ in $\rho_G$,
$j=1,\ldots , m$.
 \end{definition}
 For example, in the graph $\rho_G$ of Figure~\ref{fig:diagram1}
 we have ${\cal A}^\rho_1=\{1,5\}$ and ${\cal A}^\rho_2=\{4,7\}$. 
 The following partition summation formulas 
 extend \cite[Proposition~5.1]{LiuPrivault}
 to the counting of subgraphs with endpoints,
 and they are a special case of 
 Proposition~\ref{fjklf2} in appendix,
 which deals with joint subgraph counting. 
\begin{prop}
\label{mom-cumfor}
 Let $m\geq 0$. The moments and cumulants of $N_{y_1,\ldots , y_m}^G$ admit
 the following expressions: 
$$
  \E_\lambda \big[\big(N_{y_1,\ldots , y_m}^G\big)^n\big]=
  \sum_{\substack{\rho\in\Pi ([n]\times[r])
      \\\rho\wedge\pi=\widehat{0}} \atop {\rm (non-flat)}}
  \lambda^{|\rho |}
  \int_{(\R^d)^{|\rho|}}\prod_{\substack{ 
      1 \leq j \leq m
      \\ i\in {\cal A}^\rho_j}}
    H(x_i,y_j)
    \ \prod_{
      \substack{1 \leq k , l \le|\rho|
        \\
        \{ k , l \}\in E_{\rho_G} 
    }}H(x_k,x_l) \ \! \mu ( \mathrm{d}x_1 ) \cdots \mu ( \mathrm{d}x_{|\rho|}),
    $$
    and
    \begin{equation}
      \label{cumulant-diagram1}
    \kappa_n\big(N_{y_1,\ldots , y_m}^G\big)=
    \sum_{\substack{\rho\in\Pi_{\widehat{1}} ([n]\times[r])
        \\\rho\wedge\pi=\widehat{0}} \atop {\rm (non-flat \ \! connected)}}
  \lambda^{|\rho |}
  \int_{(\R^d)^{|\rho|}}\prod_{
    \substack{
    1 \leq j \leq m
      \\
    i\in {\cal A}^\rho_j}}  H(x_i,y_j)
  \ \prod_{\substack{
      1\leq k , l \le|\rho|
      \\
      \{ k , l \}\in E_{\rho_G} }}H(x_k,x_l) \ \! \mu ( \mathrm{d}x_1)
  \cdots \mu ( \mathrm{d}x_{|\rho|} ).
\end{equation} 
\end{prop}
We note in particular that $N_{y_1,\ldots , y_m}^G$ has positive cumulants,
and when $n=1$ the first moment of $N_{y_1,\ldots , y_m}^G$ is given by 
$$
 \E_\lambda \big[ N_{y_1,\ldots , y_m}^G \big]
 =  \lambda^r 
     \int_{(\R^d)^r}
     \prod_{
       \substack{
         1 \leq i \leq r \\
       1 \leq j \leq m
      \\
      \{v_i,w_j\} \in E_G
     }}  H(x_i,y_j)
  \ \prod_{\substack{
      1\leq k , l \leq r 
      \\
      \{ v_k , v_l \}\in E_G}}H(x_k,x_l)
  \mu ( \mathrm{d}x_1)
  \cdots
  \mu ( \mathrm{d}x_r).
$$ 
\noindent
The cumulant formula of Proposition~\ref{mom-cumfor}
is implemented in the code listed in Appendix~\ref{fjkldsf}. 
We also recall the following lemma,
see {\cite[Lemma 2.8]{LiuPrivault},
  in which maximality of connected non-flat partitions
 refers to maximizing the number of blocks.} 
\begin{lemma}
\label{numpartition}
  \noindent
  $a)$ The cardinality of the set \ $\mathcal{C} (n,r)$
 of connected non-flat partitions of $[n]\times[r]$ satisfies 
 \begin{equation}
   \label{coeff-0}
  |\mathcal{C} (n,r) | \leq n!^r r!^{n-1}, 
  \qquad n,r \geq 1. 
\end{equation}
\noindent
$b)$ 
 The cardinality of the set 
$ \mathcal{M}(n,r)$
 of maximal connected non-flat partition of $[n]\times[r]$ satisfies 
 \begin{equation}
   \nonumber 
  |\mathcal{M}(n,r)|=r^{n-1}\prod_{i=1}^{n-1}(1+(r-1)i),
  \qquad n,r\geq 1, 
\end{equation}
 with the bounds 
\begin{equation}\label{coeff-1}
    ( (r-1)r )^{n-1}(n-1)!\leq 
    |\mathcal{M}(n,r)|
     \leq ( (r-1)r )^{n-1}n!, \quad n\geq 1, \ r\geq 2. 
\end{equation}
\end{lemma}
\section{Subgraph count asymptotics} 
\label{sca}
\noindent
In this section, we let $m\geq 1$ and 
investigate the asymptotic behaviour of the cumulants
 $\kappa_n\big(N_{y_1,\ldots , y_m}^G\big)$ in \eqref{cumulant-diagram1}
 { as the intensity $\lambda$ tends to infinity,
   which extends the treatment of \cite{LiuPrivault}
   from $m=0$ to $m\geq 1$}. 
\begin{assumption}
\label{a1}
 We assume that 
\begin{enumerate}[i)]
\item $\mu$ is the Lebesgue measure on $\real^d$, and 
\item
  \label{a1-2}
  the connection function 
 $H:\R^d\times \R^d\to[0,1]$
 is translation invariant,
 i.e. $H(x,y) = H(0,y-x)$,
 $x,y \in \real^d$, and
 $$
 \int_{\real^d} H(0,y) \mathrm{d}y < \infty.
$$ 
\end{enumerate}
 \end{assumption} 
 The following result provides
growth estimates for the cumulants of $N^G_{y_1,\ldots , y_m}$. 
\begin{thm}
\label{khopone}
Let $m\geq 1$, $n\geq 1$ and $r\geq 2$, and suppose
that Assumption~\ref{a1} is satisfied.
 We have
  \begin{equation}
\label{onefixed-1}
    0<\kappa_n\big( N_{y_1,\ldots , y_m}^G \big)\leq
    n!^r r!^{n-1}
   (C \lambda )^{1+(r-1)n}
    ,
\end{equation}
and, for $n=2$, 
\begin{equation}
\label{onefixed-2}
  (r-1)r
  c^{2r}
  \lambda^{2r-1} 
  \leq \kappa_2\big(N_{y_1,\ldots , y_m}^G\big)
 \leq 
  r!
  ( C \lambda )^{2r-1}, 
\end{equation}
 where $c, C>0$ are constants independent of $r\geq 2$ and $n \geq 2$.
\end{thm}
\begin{Proof}
 According to Proposition~\ref{mom-cumfor}, every non-flat connected partition $\rho\in\Pi ([n]\times[r])$ corresponds to a summand of order $O(\lambda^{|\rho|
})$. 
 As the cardinality of {maximal}
 non-flat connected partitions is $1+(r-1)n$,
 the dominating asymptotic order is $O(\lambda^{ 1+(r-1)n })$.
 Precisely, by \eqref{coeff-0}-\eqref{coeff-1} and 
 \eqref{cumulant-diagram1}, letting
 $j_0\in \{1,\ldots , m\}$ such that ${\cal A}^\rho_{j_0} \not= \emptyset$, 
 for some $i_0 \in {\cal A}^\rho_{j_0}$ we have 
\begin{align*} 
 & 
    c^{n|E_G|}
    C^{1+(r-1)n}
        ( (r-1)r )^{n-1}(n-1)!
    \lambda^{1+(r-1)n}
    \\
        & \quad \leq      \kappa_n\big(N_{y_1,\dots,y_m}^G\big)
  \\
   & \quad
    \leq     
    \lambda^{1+(r-1)n}
    \sum_{\substack{\rho\in\Pi_{\widehat{1}} ([n]\times[r])
        \\\rho\wedge\pi=\widehat{0}} \atop {\rm (non-flat \ \! connected) \atop
            }}
    \int_{(\R^d)^{|\rho|}}\prod_{
      \substack{
        1 \leq j \leq m
      \\
    i\in {\cal A}^\rho_j}} H (x_i,y_j)
  \ \prod_{\substack{
      1\leq k , l \le|\rho|
      \\
      \{ k , l \}\in E_{\rho_G} }}H (x_k,x_l)\mathrm{d}x_1
  \cdots \mathrm{d}x_{|\rho|} 
  \\
   & \quad
    \leq     
    \lambda^{1+(r-1)n}
    \sum_{\substack{\rho\in\Pi_{\widehat{1}} ([n]\times[r])
        \\\rho\wedge\pi=\widehat{0}} \atop {\rm (non-flat \ \! connected) \atop
            }}
    \int_{(\R^d)^{|\rho|}}
     H (x_{i_0},y_{j_0})
  \ \prod_{\substack{
      1\leq k , l \le|\rho|
      \\
      \{ k , l \}\in E_{\rho_G} }}H (x_k,x_l)\mathrm{d}x_1
  \cdots \mathrm{d}x_{|\rho|} 
  \\
   & \quad
    \leq     
    \lambda^{1+(r-1)n}
    \sum_{\substack{\rho\in\Pi_{\widehat{1}} ([n]\times[r])
        \\\rho\wedge\pi=\widehat{0}} \atop {\rm (non-flat \ \! connected) \atop
            }}
    \int_{(\R^d)^{|\rho|}}
     H (x_{i_0},y_{j_0})
  \ \prod_{\substack{
      1\leq k , l \le|\rho|
      \\
      \{ k , l \}\in E_{\rho'_G} }}H (x_k,x_l)\mathrm{d}x_1
  \cdots \mathrm{d}x_{|\rho|}, 
\end{align*}
 where for every $\rho \in \Pi_{\widehat{1}} ([n]\times[r])$,
 $\rho'_G$ is a spanning tree contained in $\rho_G$,
 with vertices $\{1,\dots,|\rho|,|\rho|+j\}$
 and such that $|\rho|+j_0$ is a leaf.
 By integrating successively on the variables which
 correspond to leaves of $\rho'_G$ as in the proofs of e.g.
 Theorem~7.1 of \cite{LNS21} or
 Lemma~3.1 of \cite{can2022} and using \eqref{coeff-0}, we obtain 
\begin{align*}
\kappa_n\big(N_{y_1,\dots,y_m}^G\big)\leq (C\lambda)^{1+(r-1)n}n!^rr!^{n-1},
\end{align*}
due to Assumption~\ref{a1}-\eqref{a1-2}, 
 where $C:= \max \big( 1 , \int_{\real^d} H(0,y) dy \big)$, 
 which yields the right hand side \eqref{onefixed-1}. 
 In addition, Proposition~\ref{mom-cumfor} shows 
 that all cumulants are positive,
 which completes the proof of \eqref{onefixed-1}. 
 On the other hand, when $r\geq 2$, by
 \eqref{coeff-1} we have 
\begin{equation}
  \nonumber
  \kappa_2\big(N^G_{y_1,\ldots , y_m}\big)\geq (r-1)rC^{2r} \lambda^{2r-1},
\end{equation}
where $C>0$ is a constant independent of $r \geq 2$ and $n \geq 2$,
which shows \eqref{onefixed-2}. 
\end{Proof}
 In what follows, we consider the centered and normalized subgraph count cumulants defined as 
$$
 \widetilde{N}_{y_1,\ldots , y_m}^G
 := \frac{N_{y_1,\ldots , y_m}^G - \kappa_1 \big(N_{y_1,\ldots , y_m}^G \big)}{\sqrt{\kappa_2\big( N_{y_1,\ldots , y_m}^G \big)}}. 
$$

 \begin{corollary}
  \label{jfklds}
  Let $m\geq 1$, $n\geq 2$ and $r\geq 2$. 
 We have 
\begin{equation}
\nonumber 
\big|\kappa_n\big(\widetilde{N}_{y_1,\ldots , y_m}^G\big)\big|\leq n!^r 
C_r^{n/2}
\lambda^{-(n/2-1)},
\end{equation}
  where $C_r > 0$ is a constant depending only on $r \geq 2$.
\end{corollary}
 As a consequence of Corollary~\ref{jfklds}, 
 the skewness of $\widetilde{N}_{y_1,\ldots , y_m}^G$ satisfies  
 \begin{equation}
   \label{fjkld23} 
\big| \kappa_3 \big( \widetilde{N}_{y_1,\ldots , y_m}^G \big)\big|
\leq
C_r \lambda^{-1/2}, 
\end{equation} 
where $C_r > 0$ is a constant depending only on $r \geq 2$.
\begin{prop} 
\label{fjlfa12}
 Let $m\geq 1$. 
 The renormalized subgraph count $\widetilde{N}_{y_1,\ldots , y_m}^G$
 converges in distribution
 to the standard normal distribution
 ${\cal N}(0,1)$ as $\lambda$ tends to infinity.
\end{prop}
\begin{Proof}
  {
    { From Corollary~\ref{jfklds}
    and \eqref{fjkld23} we have 
$\kappa_1 \big( \widetilde{N}_{y_1,\ldots , y_m}^G \big) = 0$,
$\kappa_2 \big( \widetilde{N}_{y_1,\ldots , y_m}^G \big) = 1$,
and
$$
\lim_{n\to \infty}
\kappa_n \big( \widetilde{N}_{y_1,\ldots , y_m}^G \big) = 0, \qquad n\geq 3, 
$$
} 
 hence the conclusion 
 follows from Theorem~1 in \cite{Janson1988}.} 
 \end{Proof} 
 In addition, from Corollary~\ref{jfklds} and
 Lemma~\ref{Statuleviciuscond1}
 the convergence result of Proposition~\ref{fjlfa12} can be
 made more precise via the following 
 convergence bound in the Kolmogorov distance, 
{which extends Corollary~7.1 in
 \cite{LiuPrivault}
   from $m=0$ to $m\geq 1$}. 
 \begin{prop}
   \label{pkol}
  Let $m\geq 1$. We have 
\begin{equation}
\nonumber
\sup_{x\in\R}\big|\IP_\lambda \big(\widetilde{N}_{y_1,\ldots , y_m}^G \leq x\big)-\Phi(x)\big|\leq C_r \lambda^{ - 1 / ( 4r-2 )},
\qquad  r\geq 2, 
\end{equation}
where $C_r>0$ is a constant depending only on $r\geq 2$
and $\Phi$ is the cumulative distribution function of the standard normal distribution. 
 \end{prop}
  By the second moment method, see e.g. (3.4) page 54 of
\cite{jansongraphs} or Theorem 2.3.2 in \cite{roch},
we also obtain the following lower bound
for endpoint connectivity and subgraph existence. 
\begin{prop}
  \label{jklf3}
  Let $m\geq 1$.  We have 
\begin{equation}
\label{fjkl21} 
\IP_\lambda \big( 
{N}_{y_1,\ldots , y_m}^G >0 \big) \geq
\frac{\big( \E_\lambda \big[ 
    N_{y_1,\ldots , y_m}^G 
    \big]\big)^2}{ 
  \E_\lambda \big[ \big( N_{y_1,\ldots , y_m}^G \big)^2\big]},
\qquad \lambda >0. 
\end{equation}
\end{prop}
Theorem~\ref{khopone} also shows the bounds 
\begin{equation}
  \label{lb} 
\frac{C_{r,1}}{\lambda}
\leq
\frac{\kappa_2 \big(N_{y_1,\dots,y_m}^G\big)}{
    \big( \E_\lambda \big[ {N}_{y_1,\ldots , y_m}^G \big] \big)^2
    }
\leq 
\frac{C_{r,2}}{\lambda},
\qquad \lambda > 0, 
\end{equation} 
for some constants $C_{r,1},C_{r,2}>0$ depending only on $r\geq 2$,
 from which it follows that the lower bound
\eqref{fjkl21} converges to $1$ as $\lambda$ tends to infinity. 
\section{Numerical examples} 
\label{examples}
\noindent
 In this section we assume that $H$ is the Rayleigh connection function
$$
H_\beta (x,y):=e^{-\beta\|x-y\|^2}, \qquad x,y\in \real^d,
$$
where $\beta>0$, and $\mu$ is the Lebesgue measure on $\real^d$.
In this case, Assumption~\ref{a1} is satisfied.
 In the following examples, the SageMath code listed in 
 Appendix~\ref{fjkldsf} is run after loading the definitions
 of Table~\ref{t1-00}. 
  
\begin{table}[H] 
  \centering
\scriptsize 
    {
  \begin{tabular}{|ll|ll|} 
 \hline
 \multicolumn{2}{|l}{
 \EscVerb{load("cumulants_parallel.sage")}
 }
 & \multicolumn{2}{l|}{\# Loading the functions definitions ~~~~~~~~~~~~~~~~~~~~~~~~~~~~~~~~~~~~~~~~~~~~~~~~~~~~~
 }  
 \\
 \hline
 \multicolumn{2}{|l}{
 \EscVerb{λ,β = var("λ,β"); assume(β>0)}
 }
 & \multicolumn{2}{l|}{\# Variable definitions ~~~~~~~~~~~~~~~~~~~~~~~~~~~~~~~~~~~~~~~~~~~~~~~~~~~~~~~~~~~~
 }  
 \\
 \hline
 \multicolumn{2}{|l}{
 \EscVerb{def H(x,y,β): return exp(-β*(x-y)**2)}
 } 
  & \multicolumn{2}{l|}{\# Connection function $H(x,y)$}  
 \\
 \hline
 \multicolumn{2}{|l}{
 \EscVerb{def mu(x,λ,β): return 1} 
}
  & \multicolumn{2}{l|}{\# Flat intensity of $\mu (dx)$}   
 \\
\hline
\end{tabular}
}
\caption{Functions definitions.}
\label{t1-00}
\end{table} 

\vspace{-0.4cm}

\noindent
 Computations in this and the following examples
 are run on a standard desktop computer with
 an 8-core CPU at 4.10GHz. 
 {The limitations imposed by this hardware 
 configuration
 constrain the product $n\times r \times d$ to be below $15$ approximately,
 in order to maintain computation times at a reasonable level.
 The illustrations of
 figures~\ref{fig3},
     \ref{fig5-1},
     \ref{fig7} and 
     \ref{fig9}
     are provided in dimension $d=2$ for ease of visualization only.
     Actual computations may be provided in lower dimension
     $d$, due to hardware performance constraints when the cumulant
     order $n$ is beyond $4$.
  Computations in dimension $d=2$ are presented
  in Section~\ref{sub5-3} for triangles with endpoints
  and in Section~\ref{sub5-4} for trees with one endpoint
  cumulants of orders $2$ and $3$,
  while in Section~\ref{sub5-1} 
  for $3$-hop paths with two endpoint and 
  in Section~\ref{sub5-2} for four-hop paths with two endpoints
  we take $d=1$ in order to reach the cumulant
  orders $n=6$ and $n=4$. 
 In subsequent code inputs, graphs are coded by their edge set $E_G$,
 and the set of endpoints is given by the sequence 
 $\EE = [\EE_1,\ldots , \EE_m]$, where $\EE_i$ denotes the set of
 vertices of
 the subgraph $\GG$ on $\{v_1, \ldots ,v_r\}$
 which are attached to the $i$-$th$
 endpoint, $i=1,\ldots , m$, with $\EE :=[ \ ]$ the empty sequence 
 when $G$ has no endpoint ($m=0$).}  
 
\subsection{Three-hop paths with two endpoints}
\label{sub5-1}
\begin{figure}[H]
  \centering
    \includegraphics[height=4.4cm,width=0.8\linewidth]{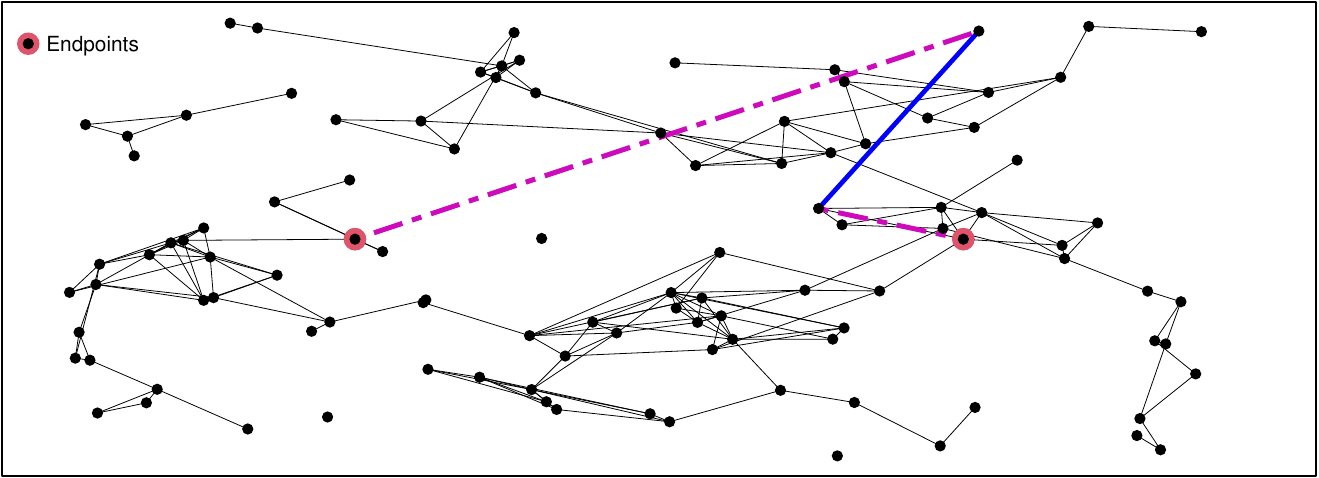} 
\caption{A $3$-hop path with two endpoints in dimension $d=2$.} 
\label{fig3}
\end{figure}

\noindent
 By a $k$-hop path, we mean a non-self intersecting path having $k$ edges. 
 We take $m=2$, $r=2$, and in Table~\ref{t1} we 
 compute the first three cumulants of $N^G_{y_1,y_2}$ when $G$ is
 a $3$-hop path with two endpoints in dimension $d=1$,
 see Figure~\ref{fig3} for an illustration in dimension $d=2$.  
 Unlike in the two-hop with two endpoints case, this
 $3$-hop count does not
 have a Poisson distribution.
 {In the following Figures~\ref{fig3},
   \ref{fig5-1},
   \ref{fig7} and
   \ref{fig9}
   the endpoints are denoted by red dots,
   and their edges are denoted by purple dashed lines. } 

\vspace{-0.3cm} 

\noindent
 To make cumulant expressions more compact,
 the exact formulas in Table~\ref{t1} 
 are expressed with $y_1=y_2=0$ and $\beta := \pi$,
 in dimension $d=1$. 

\begin{table}[H] 
  \centering
\scriptsize 
 \resizebox{\textwidth}{!}
    {
  \begin{tabular}{|ll|ll|} 
 \hline
 \multicolumn{2}{|l}{
$\GG$ = [[1,2]]; $\EE$ =[[1],[2]]; d=1 
 }
 & \multicolumn{2}{l|}{\# Single edge graph $r=2$,
 two endpoint $m=2$, dimension $d=1$~~~~~~~~~~~~~~~~~~~~~~~~~~~~~~~~~~~~~~~~~~~~~~~~~~~~~~~~~~~~~~~}   
 \\
\hline
\end{tabular}
}
  \resizebox{1\textwidth}{!}{
  \begin{tabular}{|ll|ll|} 
\hline
\multicolumn{1}{|c|}{Command} & \multicolumn{1}{c|}{Order} & \multicolumn{1}{c|}{Cumulant output} & \multicolumn{1}{c|}{Connected non-flat partitions} 
 \\ 
 \hline
\multicolumn{1}{|c|}{c(1,d,G,$\EE$,mu,H)} & \multicolumn{1}{c|}{\normalsize 1st} & \multicolumn{1}{c|}{\large $\frac{1}{\sqrt{3}} {\lambda}^{2}$} & \multicolumn{1}{c|}{\normalsize 1} 
\\
\hline
\multicolumn{1}{|c|}{c(2,d,G,$\EE$,mu,H)} & \multicolumn{1}{c|}{\normalsize 2nd} & \multicolumn{1}{c|}{\large $\left(\frac{1}{\sqrt{3}}+\frac{1}{\sqrt{2}}\right) \lambda ^3+\left(\frac{1}{\sqrt{3}}+\frac{1}{2 \sqrt{2}}\right) \lambda ^2$} & \multicolumn{1}{c|}{\normalsize 6} 
\\
\hline
\multicolumn{1}{|c|}{c(3,d,G,$\EE$,mu,H)} & \multicolumn{1}{c|}{\normalsize 3rd} & \multicolumn{1}{c|}{
  \large $\left(\sqrt{\frac{12}{7}}+\frac{3}{\sqrt{5}}+\frac{3}{\sqrt{7}}+\frac{12}{\sqrt{31}}\right) \lambda ^4+\left(\sqrt{3}+\sqrt{\frac{3}{2}}+\frac{17}{5 \sqrt{2}}+\frac{12}{\sqrt{19}}\right) \lambda^3 +
  \left(
  \frac{3}{2 \sqrt{2}}+\frac{1}{\sqrt{3}}\right)\lambda^2$} & \multicolumn{1}{c|}{\normalsize 68} 
 \\ 
\hline
\end{tabular}
}
  \caption{Cumulants of the count of $2$-hop paths with two endpoints
    in dimension $d=1$.}
\label{t1}
\end{table} 

\vspace{-0.4cm}

\noindent
Table~\ref{t1-1} lists the counts of connected non-flat partitions
and runtimes
for the computation of cumulants of orders $1$ to $6$,
and shows that such partitions represent only a fraction (around 25\%) of total partition counts.

\begin{table}[H] 
\centering
\scriptsize 
\resizebox{\textwidth}{!}{
\begin{tabular}{|c||c|c|c|c|c|c||c|c|c|} 
\hline
 Order $n$ & 2 blocks & 3 blocks & 4 blocks & 5 blocks & 6 blocks & 7 blocks & Total & $\Pi ([n]\times[2])$ & Comp. time
\\ 
\hline
 1st & 1 & 0 & 0 & 0 & 0 & 0 & 1 & 2 & 0.5s
\\ 
\hline
 2nd & 2 & 4 & 0 & 0 & 0 & 0 & 6 & 15 & 1s 
\\ 
\hline
 3rd & 4 & 32 & 32 & 0 & 0 & 0 & 68 & 203 & 3s
\\ 
\hline
 4th & 8 & 208 & 624 & 352 & 0 & 0 & 1,192 & 4,140 & 1m
\\ 
\hline
 5th & 16 & 1,280 & 8,960 & 13,904 & 5,040 & 0 & 29,200 & 115,975 & 47m
\\ 
\hline
6th & 32 & 7,744 & 116,160 & 375,776 & 351,456 & 88,544 & 939,712 & 4,213,597 & 29 hours 
\\ 
 \hline
\end{tabular}
}
\caption{Computation times and counts of connected non-flat {\em vs.} all
 partitions in $\Pi ([n]\times[2])$.}
\label{t1-1}
\end{table} 

\vspace{-0.4cm}

\begin{figure}[H]
  \centering
 \begin{subfigure}{0.49\textwidth}
   \includegraphics[width=1\linewidth, height=5cm]{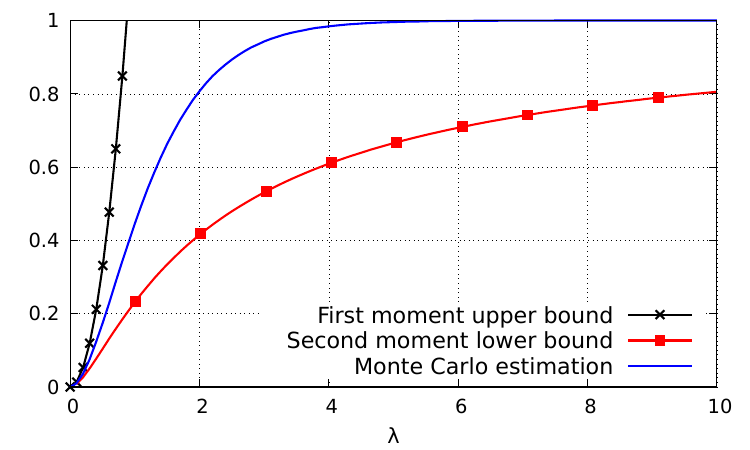} 
    \caption{First and second moment bounds \eqref{fjkl21}.} 
 \end{subfigure}
 \begin{subfigure}{0.49\textwidth}
   \vskip-0.105cm
   \includegraphics[width=1\linewidth, height=5.0cm]{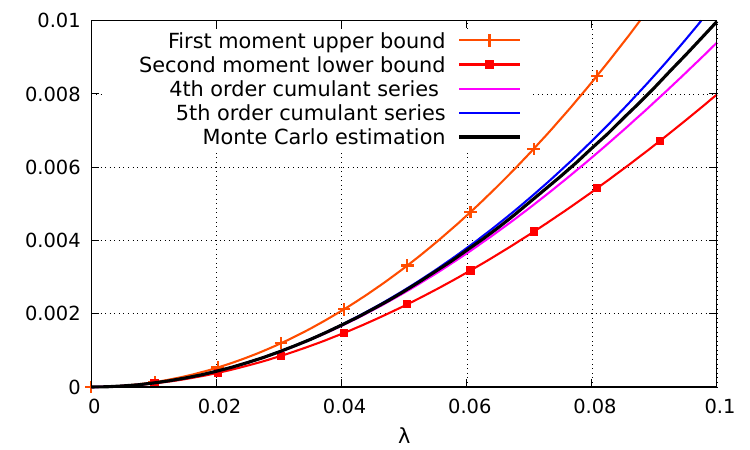} 
   \vskip-0.1cm
   \caption{Cumulant approximations \eqref{fjkl32} with $n=0$.} 
 \end{subfigure}
 \caption{Connection probabilities.} 
\label{fig2-111} 
\end{figure}

\noindent
 Figure~\ref{fig2-111} presents connectivity estimates based
 on the moment 
 and cumulant formulas of Propositions~\ref{jklf3} and \ref{fdshkf0},
 in dimension $d=1$. 

\noindent
\subsection{Four-hop paths with two endpoints}
\label{sub5-2}
\begin{figure}[H]
\centering
\includegraphics[height=4.4cm,width=0.8\linewidth]{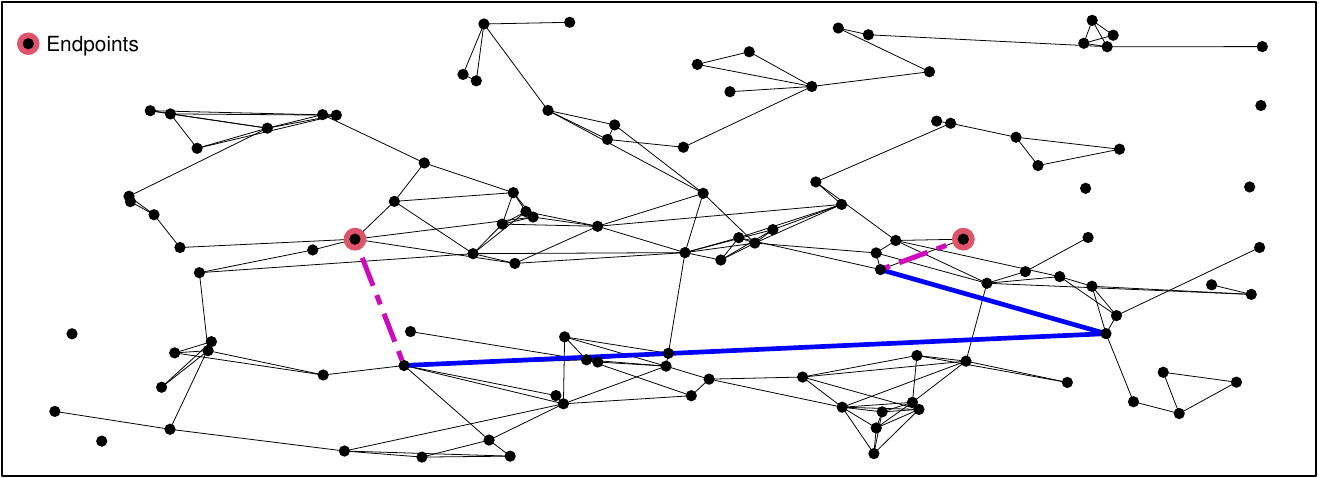} 
\caption{A $4$-hop path with two endpoints in dimension $d=2$.} 
\label{fig5-1}
 \end{figure}

\noindent 
 Here, we take $m=2$ and $r=3$, and in Table~\ref{t2} we 
 compute the first cumulant of $N^G_{y_1,y_2}$ when $G$ is
 a four-hop path with two endpoints in dimension $d=1$, 
 see Figure~\ref{fig5-1} for an illustration in dimension $d=2$.  

\vspace{-0.3cm} 

\noindent
 The closed-form expressions in Table~\ref{t2} 
 are expressed with $y_1=y_2=0$ and $\beta := \pi$,
 in dimension $d=1$. 

 \begin{table}[H] 
  \centering
\resizebox{\textwidth}{!}
      {
  \begin{tabular}{|ll|ll|} 
 \hline
 \multicolumn{2}{|l}{
\large $\GG$ = [[1,2],[2,3]];  $\EE$=[[1],[3]]; d=1
}
 & \multicolumn{2}{l|}{\large \# $4$-hop path $r=3$,
 two endpoints $m=2$, dimension $d=1$}   
 \\
\hline
\hline
\multicolumn{1}{|c|}{Instruction} & \multicolumn{1}{c|}{Order} & \multicolumn{1}{c|}{Cumulant output} & \multicolumn{1}{c|}{Connected non-flat partitions} 
 \\ 
 \hline
\multicolumn{1}{|c|}{c(1,d,G,$\EE$,mu,H)} & \multicolumn{1}{c|}{1st} & \multicolumn{1}{c|}{\Large $\frac{{\lambda}^{3}}{2}$} & \multicolumn{1}{c|}{\large 1} 
 \\ 
 \hline
\multicolumn{1}{|c|}{c(2,d,G,$\EE$,mu,H)} & \multicolumn{1}{c|}{2nd} & \multicolumn{1}{c|}{\Large $\frac{1}{6} \lambda ^3 \left(\left(\sqrt{6}+4 \sqrt{\frac{3}{5}}+\frac{3}{2 \sqrt{2}}+\frac{12}{\sqrt{7}}\right) \lambda ^2+\left(3 \sqrt{3}+16 \sqrt{\frac{3}{7}}+8 \sqrt{\frac{3}{11}}+\frac{3}{2 \sqrt{2}}+\frac{6}{\sqrt{5}}\right) \lambda +\sqrt{3}+\sqrt{6}+6\right)$}
        & \multicolumn{1}{c|}{\large 33} 
\\ 
\hline
\end{tabular}
}
\caption{First and second cumulants of the count of four-hop paths with two endpoints.}
\label{t2} 
\end{table} 

\vspace{-0.4cm}

\noindent
Table~\ref{t1-1-2} presents the counts of connected non-flat partitions at different orders, and shows that such partitions represent only a fraction (around 10\%) of total partition counts. 
 
\begin{table}[H] 
  \centering
\scriptsize 
\resizebox{\textwidth}{!}{
\begin{tabular}{|c||c|c|c|c|c|c|c||c|c|c|} 
\hline
 Order $n$ & 3 blocks & 4 blocks & 5 blocks & 6 blocks & 7 blocks & 8 blocks & 9 blocks & Total & $\Pi ([n]\times[3])$ & Comp. time
\\ 
\hline
 1st & 1 & 0 & 0 & 0 & 0 & 0 & 0 & 1 & 5 & 1s
\\ 
\hline
 2nd & 6 & 18 & 9 & 0 & 0 & 0 & 0 & 33 & 203 & 2s 
\\ 
\hline
 3rd & 36 & 540 & 1242 & 864 & 189 & 0 & 0 & 2871 & 21,147 & 4m 
\\ 
\hline
 4th & 216 & 13,608 & 94,284 & 186,624 & 145,908 & 48,276 & 5,589 & 494,500 & 4,213,597 & 19 hours 
\\ 
 \hline
\end{tabular}
}
\caption{Computation times and counts of connected non-flat {\em vs.} all partitions in $\Pi ([n]\times[3])$.}
\label{t1-1-2}
\end{table} 

\vspace{-0.4cm}

\begin{figure}[H]
  \centering
 \begin{subfigure}[b]{0.49\textwidth}
    \includegraphics[width=1\linewidth, height=4.8cm]{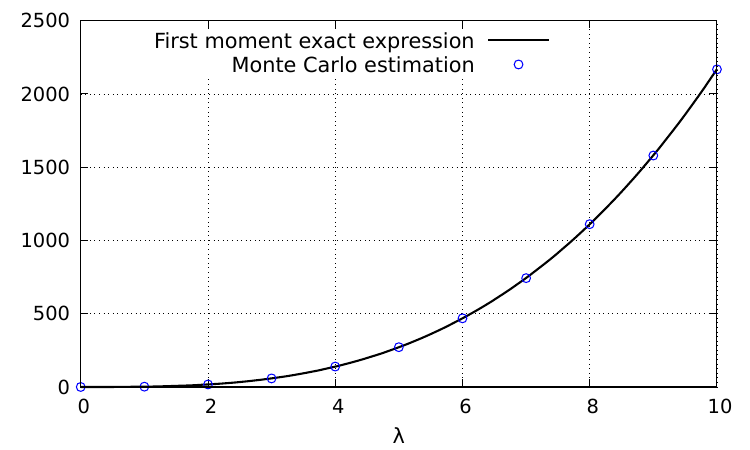} 
    \caption{First moment.} 
 \end{subfigure}
 \begin{subfigure}[b]{0.49\textwidth}
    \includegraphics[width=1\linewidth, height=4.8cm]{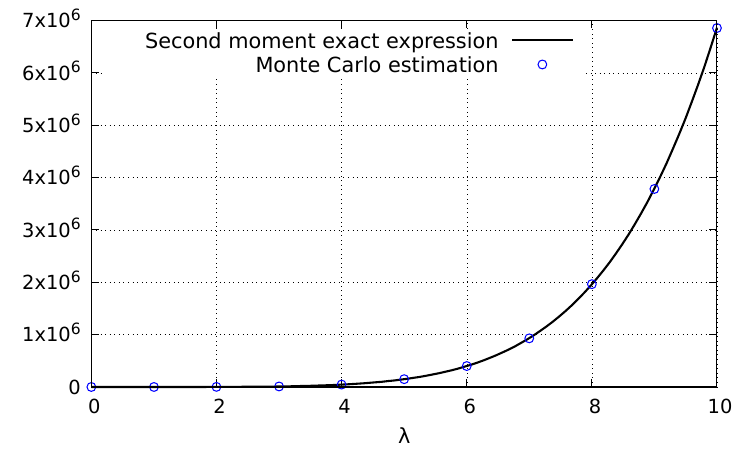} 
    \caption{Second moment.} 
 \end{subfigure}
 \begin{subfigure}[b]{0.50\textwidth}
 \hskip-0.5cm
    \includegraphics[width=1.06\linewidth, height=4.8cm]{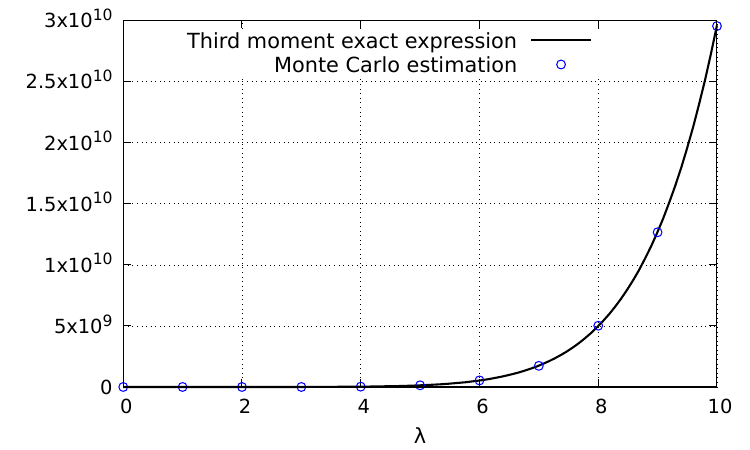} 
    \caption{Third moment.} 
 \end{subfigure}
  \centering
  \begin{subfigure}[b]{0.49\textwidth}
  \hskip-0.3cm
    \includegraphics[width=1.03\linewidth, height=4.8cm]{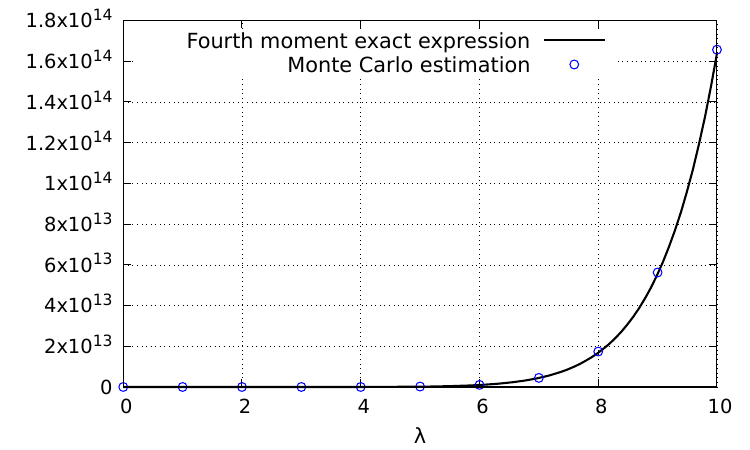} 
    \caption{Fourth moment.} 
 \end{subfigure}
  \caption{Moment estimates.} 
\label{fig2-11} 
\end{figure}

\noindent
 In Figure~\ref{fig2-11} we plot the corresponding
 moment expressions {\em vs.} their Monte Carlo estimates
 in dimension $d=1$, {with the parameters of Table~\ref{t2}}.
 
\subsection{Triangles with endpoints}
\label{sub5-3}
\begin{figure}[H]
\centering
\includegraphics[width=0.8\linewidth]{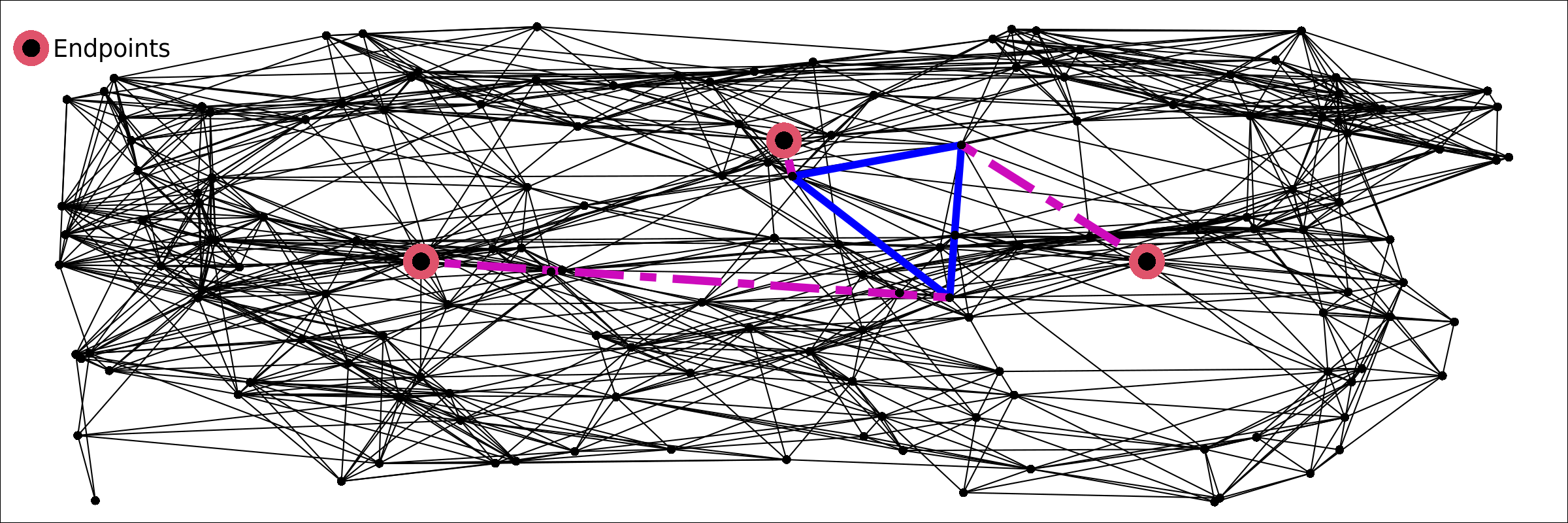} 
\caption{A triangle with three endpoints in dimension $d=2$.} 
\label{fig7}
\end{figure}

\noindent
 Taking $r=3$ and $m=3$, in Table~\ref{t3} we 
 compute the first and second
 cumulants of $N^G_{y_1,y_1,y_3}$ when $G$ is
 a triangle with three endpoints in dimension $d=1$, 
 see Figure~\ref{fig7} for an illustration in dimension $d=2$.  

\vspace{-0.3cm} 

\noindent
 The closed-form expressions in Table~\ref{t3} 
 are expressed with $y_1=y_2=y_3=0$ and $\beta := \pi$,
 in dimension $d=1$. 

 \begin{table}[H] 
  \centering
    \resizebox{\textwidth}{!}
      {
  \begin{tabular}{|ll|ll|} 
 \hline
 \multicolumn{2}{|l}{
   \large $\GG$ = [[1,2],[2,3],[3,1]];
    $\EE$=[[1],[2],[3]]; d=1;
}
  & \multicolumn{2}{l|}{\large \# Triangle graph $r=3$; three endpoints $m=3$; dimension $d=1$} 
 \\
\hline
\hline
\multicolumn{1}{|c|}{Instruction} & \multicolumn{1}{c|}{Order} & \multicolumn{1}{c|}{Cumulant output} & \multicolumn{1}{c|}{Connected non-flat partitions} 
\\ 
\hline
\multicolumn{1}{|c|}{c(1,d,G,$\EE$,mu,H)} & \multicolumn{1}{c|}{1st} & \multicolumn{1}{c|}{\Large $\frac{{\lambda}^{3}}{4}$} & \multicolumn{1}{c|}{1} 
\\ 
\hline
\multicolumn{1}{|c|}{c(2,d,G,$\EE$,mu,H)} & \multicolumn{1}{c|}{2nd} & \multicolumn{1}{c|}{\Large ~~~~$
  \left(
  \frac{\sqrt{3}}{8}   + \frac{3}{8}
  \right)
        {\lambda}^{5}
        +
        \left(
        \frac{2\sqrt{105} }{35}
  + \frac{\sqrt{3} }{5}
  + \frac{3}{4}
  \right) {\lambda}^{4} 
  +
  \left(
  \frac{3\sqrt{35} }{35}
  + \frac{\sqrt{2} }{5}
  + \frac{1}{4}
  \right) {\lambda}^{3}
  $~~~~} & \multicolumn{1}{c|}{33} 
\\ 
\hline
\end{tabular}
}
\caption{First and second cumulants of the count of triangles with three endpoints.}
\label{t3}
\end{table} 

\vspace{-0.4cm}

\noindent
Table~\ref{t1-1-3} presents computation times in dimension $d=2$.
 
\vspace{-0.1cm}

\begin{table}[H] 
  \centering
\scriptsize 
\resizebox{\textwidth}{!}{
\begin{tabular}{|c||c|c|c|c|c||c|c|c|} 
\hline
 Order $n$ & 3 blocks & 4 blocks & 5 blocks & 6 blocks & 7 blocks & Total & $\Pi ([n]\times[3])$ & Comp. time
\\ 
\hline
 1st & 1 & 0 & 0 & 0 & 0 & 1 & 5 & 1s
\\ 
\hline
 2nd & 6 & 18 & 9 & 0 & 0 & 33 & 203 & 21s
\\ 
\hline
 3rd & 36 & 540 & 1,242 & 864 & 189 & 2,871 & 21,147 & 1 hour 
\\ 
 \hline
\end{tabular}
}
\caption{Computation times and counts of connected non-flat {\em vs.} all partitions in $\Pi ([n]\times[3])$.}
\label{t1-1-3}
\end{table} 

\vspace{-0.4cm}

\begin{figure}[H]
\centering
\begin{subfigure}{.5\textwidth}
\centering
\includegraphics[width=1.\textwidth]{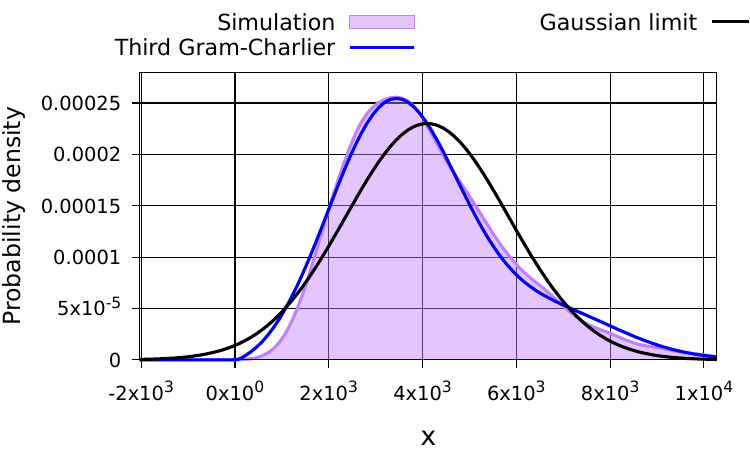} 
\vskip-0.08cm
\caption{\small $\lambda = 50$.} 
\end{subfigure}
\hskip-0.2cm
\begin{subfigure}{.5\textwidth}
\centering
\includegraphics[width=1.\textwidth]{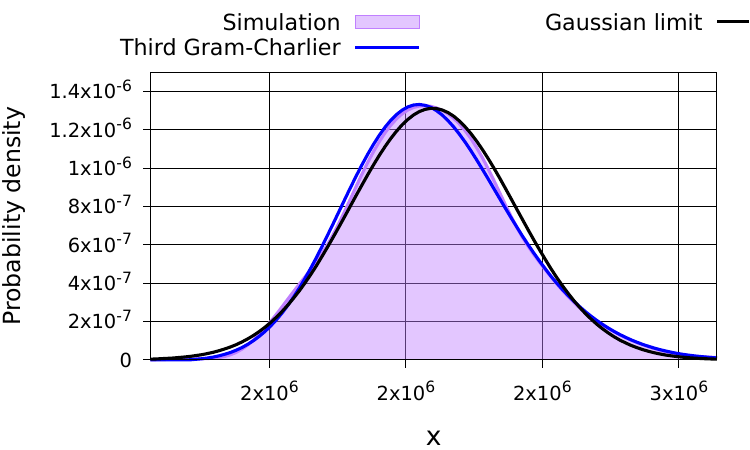} 
\vskip-0.08cm
\caption{\small $\lambda = 400$.} 
\end{subfigure}
\vskip-0.16cm
\caption{\small Gram-Charlier density expansions {\em vs.} Monte Carlo density estimation.} 
\label{fig5}
\end{figure}

\noindent
 Figure~\ref{fig5} presents second and third order Gram-Charlier
 expansions \eqref{gram_charlier-1}-\eqref{gram_charlier-2} 
 for the probability density function of the
 count $N^G_{y_1,y_1,y_3}$ of triangles with three endpoints,
 based on exact second and third cumulant expressions. 
 
\vspace{-0.3cm}

\noindent
In Figure~\ref{fig5}, the purple areas correspond to probability density estimates obtained by Monte Carlo simulations with $\beta = 2$. 
The second order expansions correspond to
the Gaussian diffusion approximation 
obtained by matching first and second order moments. 
Figure~\ref{fig5} shows that
the actual probability density estimates obtained by simulation 
can be significantly different from
their Gaussian diffusion approximations when 
skewness takes large absolute values. 
In addition, in Figure~\ref{fig5} 
the fourth order Gram-Charlier expansions appear to give the best fit
to the actual probability densities, 
which have positive skewness. 

\subsection{Trees with one endpoint}
\label{sub5-4}
\begin{figure}[H]
\centering
\includegraphics[height=4.4cm,width=0.8\linewidth]{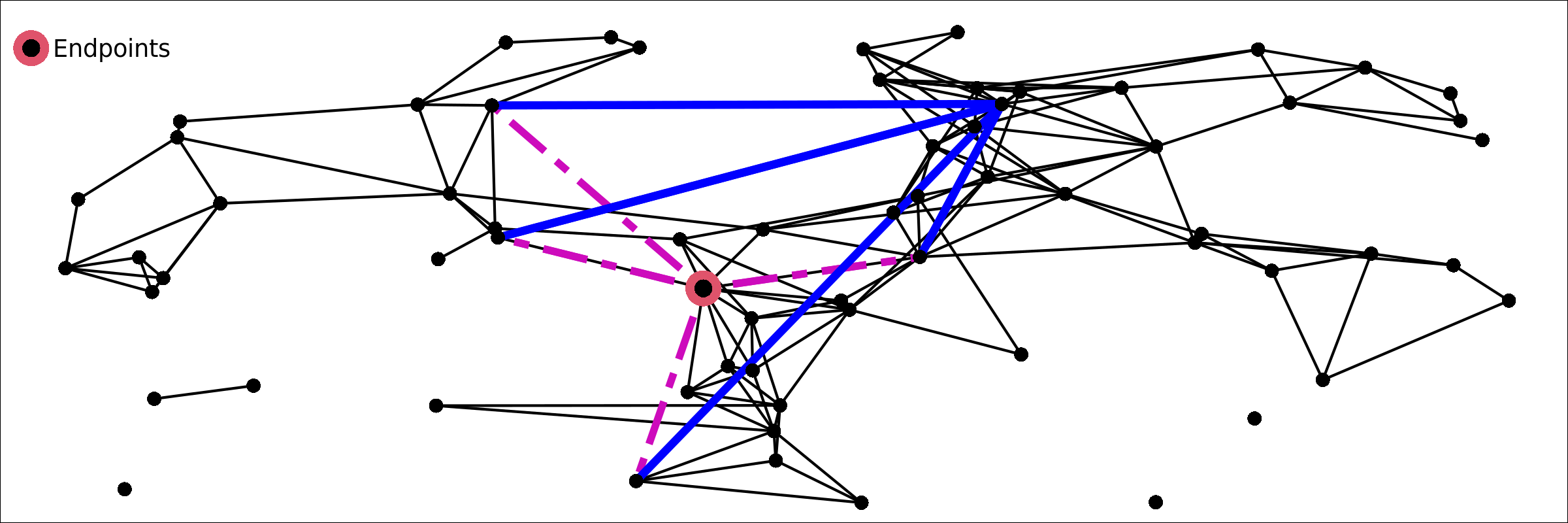} 
\caption{Four trees with a single endpoint in dimension $d=2$.} 
\label{fig9}
\end{figure}

\noindent
 Here we take $r=4$ and $m=1$, and in Table~\ref{t3-1} we 
 compute the first and second cumulants of $N^G_{y_1,y_1,y_3}$ when $G$ is
 made of a tree and a single endpoint in dimension $d=2$, 
 see Figure~\ref{fig9} for an illustration. 
    
\vspace{-0.3cm} 

\noindent
{
In Table~\ref{t3-1} 
  the location of the unique endpoint has no
  impact on cumulant expressions 
  due to space homogeneity of the underlying
  Poisson point process.} 

    \begin{table}[H] 
  \centering
    \resizebox{\textwidth}{!}
      {
  \begin{tabular}{|ll|ll|} 
 \hline
 \multicolumn{2}{|l}{
   $\GG$ = [[1,2],[2,3],[2,4]];
    $\EE$=[[1,3,4]]; d=2;
}
  & \multicolumn{2}{l|}{\# Tree $r=4$; single endpoint $m=1$; dimension $d=2$} 
 \\
\hline
\hline
\multicolumn{1}{|c|}{Instruction} & \multicolumn{1}{c|}{Order} & \multicolumn{1}{c|}{Cumulant output} & \multicolumn{1}{c|}{Connected non-flat partitions} 
\\ 
\hline
\multicolumn{1}{|c|}{c(1,d,G,$\EE$,mu,H)} & \multicolumn{1}{c|}{\small 1st} & \multicolumn{1}{c|}{\large $\frac{{\lambda}^{4}}{12}$} & \multicolumn{1}{c|}{\small 1} 
\\ 
\hline
\multicolumn{1}{|c|}{c(2,d,G,$\EE$,mu,H)} & \multicolumn{1}{c|}{\small 2nd} & \multicolumn{1}{c|}{\large ~~~~$
\frac{41{\lambda}^{7}}{384}  + \frac{99039{\lambda}^{6}}{165760}  + \frac{232885{\lambda}^{5}}{175824}  + \frac{37{\lambda}^{4}}{50}$~~~~} & \multicolumn{1}{c|}{\small 208} 
\\ 
\hline
\end{tabular}
}
\caption{First and second cumulants of the count of trees with one endpoint.}
\label{t3-1}
\end{table} 

\vspace{-0.4cm}

\noindent
The computation times presented in 
Table~\ref{t1-1-4} are for dimension $d=2$.
 
\vspace{-0.2cm}

\begin{table}[H] 
  \centering
\scriptsize 
\resizebox{\textwidth}{!}{
\begin{tabular}{|c||c|c|c|c|c|c|c||c|c|c|} 
\hline
 Order $n$ & 4 blocks & 5 blocks & 6 blocks & 7 blocks & 8 blocks & 9 blocks & 10 blocks & Total & $\Pi ([n]\times[4])$ & Comp. time
\\ 
\hline
 1st & 1 & 0 & 0 & 0 & 0 & 0 & 0 & 1 & 15 & 1s
\\ 
\hline
 2nd & 24 & 96 & 72 & 15 & 0 & 0 & 0 & 208 & 4,140 & 2m
\\ 
\hline
 3rd & 576 & 13,824 & 50,688 & 59,904 & 29,952 & 6,912 & 640 & 162,496 & 4,213,597 & 40 hours
\\ 
 \hline
\end{tabular}
}
\caption{Computation times and counts of connected non-flat {\em vs.} all partitions in $\Pi ([n]\times[4])$.} 
\label{t1-1-4}
\end{table} 

\vspace{-0.4cm}

\begin{figure}[H]
  \centering
 \begin{subfigure}[b]{0.50\textwidth}
 \hskip-0.5cm
    \includegraphics[width=1.06\linewidth, height=5cm]{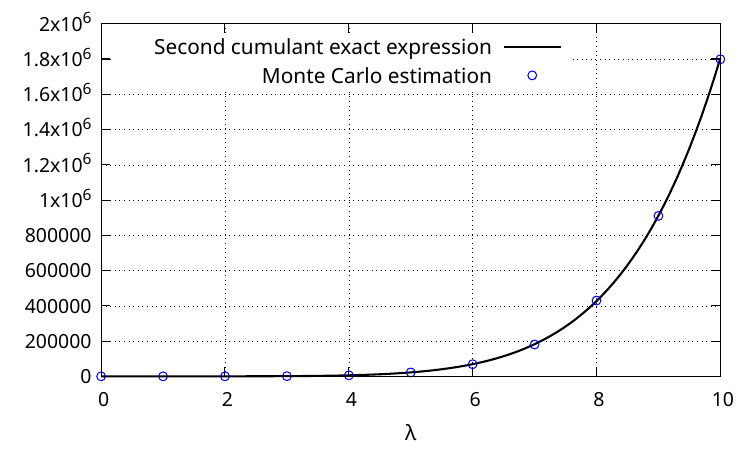} 
    \caption{Second cumulant.} 
 \end{subfigure}
  \centering
  \begin{subfigure}[b]{0.49\textwidth}
  \hskip-0.3cm
    \includegraphics[width=1.03\linewidth, height=5cm]{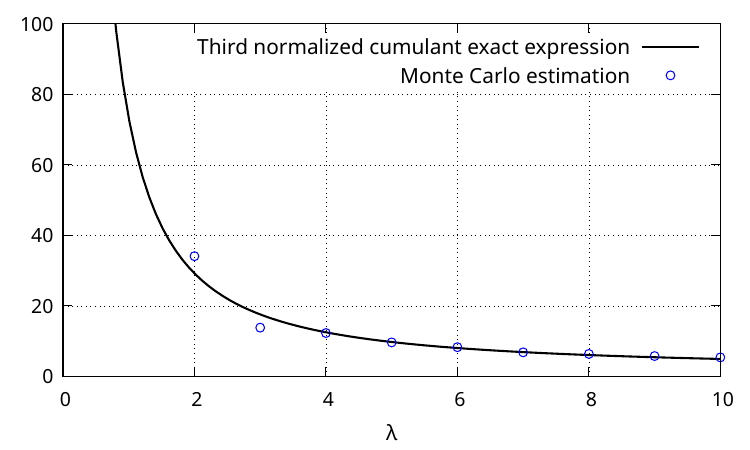} 
    \caption{Normalized third cumulant.} 
 \end{subfigure}
\caption{Cumulant estimates.} 
\label{fig2-11-2} 
\end{figure}

\noindent
In Figure~\ref{fig2-11-2} we plot the
second cumulant of ${N}_{y_1,\ldots , y_m}^G$ 
and the third cumulant of  $\widetilde{N}_{y_1,\ldots , y_m}^G$
 {\em vs.} their Monte Carlo estimates
 in dimension $d=2$ with $y_1 = y_2 = y_3 = 0$. 

\subsection{Correlation of triangles {\em vs.} four-hop counts} 
\label{corrof}
\noindent
In this example, we run the joint cumulant code provided
in Appendix~\ref{fjkldsf-2}
to compute the correlation of triangle and four-hop counts
without endpoints, as a function of the intensity parameter $\lambda$.
 Here, $\mu$ is taken to be a finite measure
 as no endpoints are considered, i.e.
 we have $\EE$=[ \hskip0.03cm ] and $m=0$,
 and the SageMath code listed in 
 Appendix~\ref{fjkldsf-2} is run after loading the definitions
 of Table~\ref{t1-002}. 
 
\begin{table}[H] 
  \centering
\scriptsize 
 \resizebox{\textwidth}{!}
    {
  \begin{tabular}{|ll|ll|} 
 \hline
 \multicolumn{2}{|l}{
 \EscVerb{load("cumulants_parallel.sage");load("jointcumulants.sage")}
 }
 & \multicolumn{2}{l|}{\# Loading the functions definitions ~~~~~~~~~~
 }  
 \\
 \hline
 \multicolumn{2}{|l}{
 \EscVerb{λ,β = var("λ,β"); assume(β>0)}
 }
 & \multicolumn{2}{l|}{\# Variable definitions ~~~~~~~~~~~~~~~~~~~
 }  
 \\
 \hline
 \multicolumn{2}{|l}{
 \EscVerb{def H(x,y,β): return exp(-β*(x-y)**2)}
 } 
  & \multicolumn{2}{l|}{\# Connection function $H(x,y)$}  
 \\
 \hline
 \multicolumn{2}{|l}{
 \EscVerb{def mu(x,λ,β): return exp(-β*x**2)} 
}
  & \multicolumn{2}{l|}{\# Finite intensity measure $\mu (dx)$}   
 \\
\hline
\end{tabular}
}
\caption{Functions definitions.}
\label{t1-002}
\end{table} 

\vspace{-0.4cm}

\noindent
 The closed-form expressions in Table~\ref{t3-1-1}
 are expressed with $\beta := \pi$, in dimension $d=2$. 

\begin{table}[H] 
  \centering
    \resizebox{1.0\textwidth}{!}
      {
  \begin{tabular}{|ll|ll|} 
 \hline
 \multicolumn{4}{|l|}{
$\GG 1$ = [[1,2],[2,3],[3,1]]; 
$\GG 2$ = [[1,2],[2,3],[3,4],[4,5]]; 
$\GG 2$c = [[4,5],[5,6],[6,7],[7,8]]; 
$\GG$ = [$\GG 1$,$\GG 2$c]; $\EE$=[]; d=2;
}
 \\
 \hline
 \multicolumn{4}{|l|}{\# Triangles G1 and $4$-hops G2;
    $r_1=3$, $r_2=5$; no endpoints $m=0$; dimension $d=2$} 
 \\
\hline
\hline
\multicolumn{1}{|c|}{Instruction} & \multicolumn{1}{c|}{Order} & \multicolumn{1}{c|}{Cumulant output} & \multicolumn{1}{c|}{Connected non-flat partitions} 
\\ 
\hline
\multicolumn{1}{|c|}{c(2,d,G1,$\EE$,mu,H)} & \multicolumn{1}{c|}{\small 2nd} & \multicolumn{1}{c|}{\normalsize \Large$\frac{3 \lambda^5}{64} + \frac{6\lambda^4}{25} + \frac{3\lambda^3}{8}$} & \multicolumn{1}{c|}{\small 33} 
\\ 
\hline
\multicolumn{1}{|c|}{c(2,d,G2,$\EE$,mu,H)} & \multicolumn{1}{c|}{\small 2nd} & \multicolumn{1}{c|}{\Large$
  \frac{7344738590701\lambda^9}{687218605505250} + \cdots 
    $} & \multicolumn{1}{c|}{\small 1545} 
\\ 
\hline
\multicolumn{1}{|c|}{jc(d,G,$\EE$,mu,H)} & \multicolumn{1}{c|}{\small 2nd joint} & \multicolumn{1}{c|}{\Large ~~~~~~$
  \frac{34409 \lambda^7}{1537920} +
  \frac{9101145477 \lambda^6}{55004486680}
  + \frac{10774977 \lambda^5}{28148120}
  $~~~~~} & \multicolumn{1}{c|}{\small 135} 
\\
\hline
\end{tabular}
}
\caption{Second (joint) moments of triangle counts {\em vs.} four-hop counts}\label{t3-1-1}
\end{table} 

\vspace{-0.4cm}

\begin{figure}[H]
  \centering
 \begin{subfigure}[b]{0.50\textwidth}
 \hskip-0.5cm
    \includegraphics[width=1.06\linewidth, height=5cm]{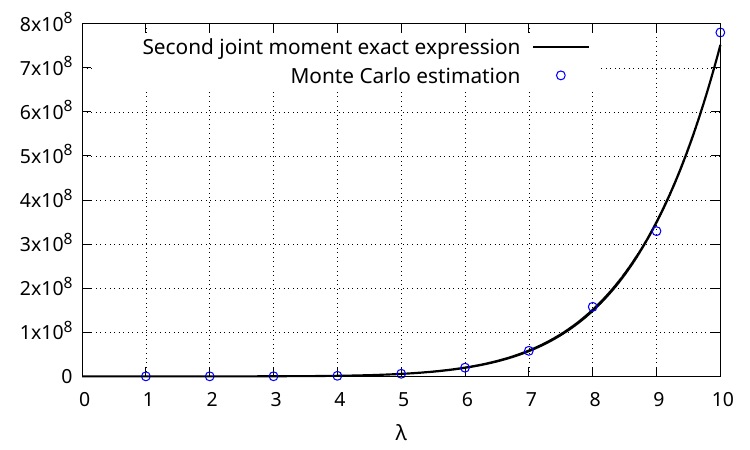} 
    \caption{Second joint cumulant.} 
 \end{subfigure}
  \centering
  \begin{subfigure}[b]{0.49\textwidth}
  \hskip-0.3cm
    \includegraphics[width=1.03\linewidth, height=5cm]{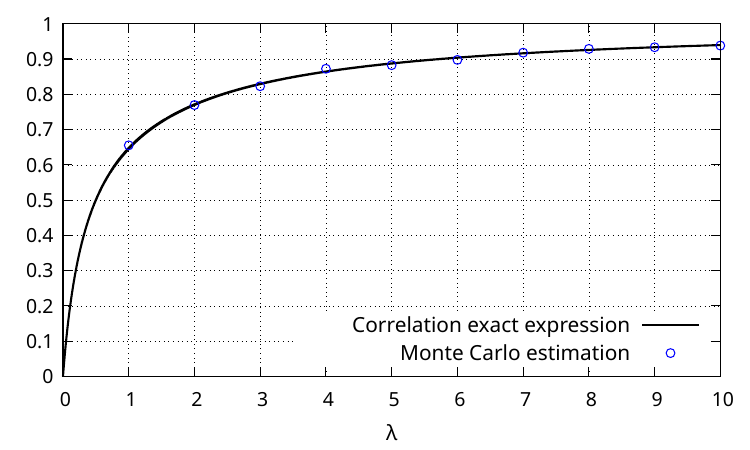} 
    \caption{Correlation.} 
 \end{subfigure}
\caption{Correlation and second joint cumulant estimates.} 
\label{fig2-11-3} 
\end{figure}

\noindent
In Figure~\ref{fig2-11-3} we plot the 
second joint cumulant and correlation of $\big({N}^{G_1},{N}^{G_2}\big)$ 
 {\em vs.} their Monte Carlo estimates
 in dimension $d=1$. 

\noindent
The limit correlation as $\lambda$ tends to infinity can be
exactly estimated from Table~\ref{t3-1-1} as 
$$
\frac{34409 }{1537920} \sqrt{\frac{64}{3} \times
  \frac{687218605505250}{7344738590701}} \approx 0.999602.
$$ 

\appendix

\section{Multivariate moment and cumulant formulae}
\label{appendixa}
\noindent 
In this section, we prove an extension of Proposition~\ref{mom-cumfor}
for the joint moments and cumulants of subgraph counts. 
The next definition extends Definition~\ref{def-1}.
\begin{definition}
Given $r_1,\dots,r_n \geq 1$, we set 
$$
  \pi_i=\left\{(i,1), \ldots ,(i,r_i)\right\},
  \quad
  i=1, \ldots , n, 
$$
 and $\pi := \{ \pi_1,\ldots , \pi_n \}$. 
\begin{enumerate}[i)]
   \item A set partition $\sigma\in\Pi( \pi_1 \cup \cdots \cup \pi_n )$ is connected if $\sigma\vee\pi=\widehat{1}$.      
\item 
 A set partition $\sigma\in\Pi ( \pi_1 \cup \cdots \cup \pi_n )$ is non-flat if $\sigma\wedge\pi=\widehat{0}$. 
\end{enumerate} 
\noindent
 We let $\Pi_{\widehat{1}}( \pi_1 \cup \cdots \cup \pi_n )$ denote the collection of all connected partitions of $\pi_1 \cup \cdots \cup \pi_n$. 
\end{definition}
  In what follows, every partition
  $\rho \in \Pi(\pi_1\cup \cdots \cup \pi_n )$
  will be arranged into a 
 diagram denoted by $\Gamma(\rho ,\pi)$, 
 by arranging $\pi_1,\dots,\pi_n$ into $n$ rows 
 and connecting together the elements of every block of $\rho$, 
 see 
 Figure~\ref{fig:diagram2} for two illustrations with
 $n=5$, 
 $(r_1, r_2 , r_3 , r_4 , r_5) = (3,2,4,3,4)$. 
 
\begin{figure}[H]
\captionsetup[subfigure]{font=footnotesize}
\centering
\subcaptionbox{Non-connected partition diagram $\Gamma(\rho,\pi)$.}[.5\textwidth]{
\begin{tikzpicture}[scale=0.8] 
\draw[black, thick] (0,0) rectangle (5,6);
\node[anchor=east,font=\small] at (0.8,5) {1};
\node[anchor=east,font=\small] at (0.8,4) {2};
\node[anchor=east,font=\small] at (0.8,3) {3};
\node[anchor=east,font=\small] at (0.8,2) {4};
\node[anchor=east,font=\small] at (0.8,1) {5};
\node[anchor=south,font=\small] at (1,0) {1};
\node[anchor=south,font=\small] at (2,0) {2};
\node[anchor=south,font=\small] at (3,0) {3};
\node[anchor=south,font=\small] at (4,0) {4};
\filldraw [gray] (1,1) circle (2pt);
\filldraw [gray] (2,1) circle (2pt);
\filldraw [gray] (3,1) circle (2pt);
\filldraw [gray] (4,1) circle (2pt);
\filldraw [gray] (1,2) circle (2pt);
\filldraw [gray] (2,2) circle (2pt);
\filldraw [gray] (3,2) circle (2pt);
\filldraw [gray] (1,3) circle (2pt);
\filldraw [gray] (2,3) circle (2pt);
\filldraw [gray] (3,3) circle (2pt);
\filldraw [gray] (4,3) circle (2pt);
\filldraw [gray] (2,3) circle (2pt);
\filldraw [gray] (1,4) circle (2pt);
\filldraw [gray] (2,4) circle (2pt);
\filldraw [gray] (1,5) circle (2pt);
\filldraw [gray] (2,5) circle (2pt);
\filldraw [gray] (3,5) circle (2pt);
\draw[very thick] (1,5) -- (1,4);
\draw[very thick] (2,5) -- (2,4);
\draw[very thick] (2,5) -- (3,5);
\draw[very thick] (1,2) -- (1,1);
\draw[very thick] (2,3) -- (2,2);
\draw[very thick] (2,1) -- (3,1) -- (4,1);
\draw[very thick] (3,2) -- (4,3);
 \begin{pgfonlayer}{background}
    \filldraw [line width=4mm,black!3]
      (0.27,0.27)  rectangle (4.76,5.76);
  \end{pgfonlayer}
\end{tikzpicture}}
\subcaptionbox{Connected partition diagram $\Gamma(\rho,\pi)$.}[.49\textwidth]{
\begin{tikzpicture}[scale=0.8] 
\draw[black, thick] (0,0) rectangle (5,6);

\node[anchor=east,font=\small] at (0.8,5) {1};
\node[anchor=east,font=\small] at (0.8,4) {2};
\node[anchor=east,font=\small] at (0.8,3) {3};
\node[anchor=east,font=\small] at (0.8,2) {4};
\node[anchor=east,font=\small] at (0.8,1) {5};

\node[anchor=south,font=\small] at (1,0) {1};
\node[anchor=south,font=\small] at (2,0) {2};
\node[anchor=south,font=\small] at (3,0) {3};
\node[anchor=south,font=\small] at (4,0) {4};

\filldraw [gray] (1,1) circle (2pt);
\filldraw [gray] (2,1) circle (2pt);
\filldraw [gray] (3,1) circle (2pt);
\filldraw [gray] (4,1) circle (2pt);
\filldraw [gray] (1,2) circle (2pt);
\filldraw [gray] (2,2) circle (2pt);
\filldraw [gray] (3,2) circle (2pt);
\filldraw [gray] (1,3) circle (2pt);
\filldraw [gray] (2,3) circle (2pt);
\filldraw [gray] (3,3) circle (2pt);
\filldraw [gray] (4,3) circle (2pt);
\filldraw [gray] (2,3) circle (2pt);
\filldraw [gray] (1,4) circle (2pt);
\filldraw [gray] (2,4) circle (2pt);
\filldraw [gray] (1,5) circle (2pt);
\filldraw [gray] (2,5) circle (2pt);
\filldraw [gray] (3,5) circle (2pt);

\draw[very thick] (1,5) -- (1,4); 

\draw[very thick] (1,2) -- (1,1);
\draw[very thick] (2,2) -- (2,4);
\draw[very thick] (2,1) -- (3,2) -- (4,3);

 \begin{pgfonlayer}{background}
    \filldraw [line width=4mm,black!3]
      (0.27,0.27)  rectangle (4.76,5.76);
  \end{pgfonlayer}
\end{tikzpicture}}%
\caption{Two examples of partition diagrams.}
\label{fig:diagram2}
\end{figure}
\vspace{-0.4cm}

\noindent
 Definition~\ref{part-1} extends \cite[Definition~2.4]{LiuPrivault}
to the multivariate setting. 
\begin{definition}\label{part-1}
    ~~
  \begin{enumerate}[\rm 1)]
\item
  Given $\rho\in\Pi(\pi_1\cup \cdots \cup \pi_n )$,
  we let $\sigma_\rho$ be the partition of $[n]$ defined by
  the condition 
  $$\rho\vee\pi=\bigg\{\bigcup_{i\in b}\pi_i:b\in\sigma_\rho \bigg\}. 
$$
\item For any non-empty set $b \subset [n]$, we let 
$$\rho_b:=\bigg\{c\in\rho:c\subset \bigcup_{i\in b}\pi_i\bigg\}.
$$
\end{enumerate}
\end{definition}
 As an example, in Figure~\ref{fig:diagram3}-$a)$, when $b = \{1,2\}$ we have
$$
\rho_{\{1,2\}} = \big\{\{(1,1),(2,1)\}, \{(1,2),(1,3),(2,2)\}\big\}. 
$$

\tikzset{hide labels/.style={every label/.append style={text opacity=0}}}

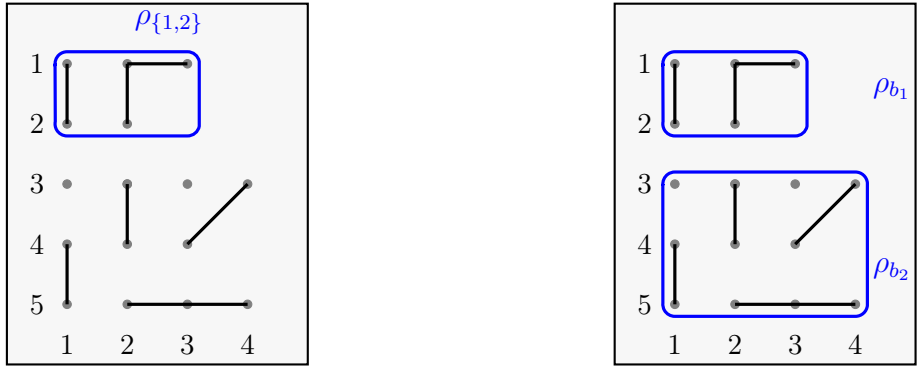
\begin{figure}[H]
\captionsetup[subfigure]{font=footnotesize}
\centering
\subcaptionbox{Connected subpartition $\rho_{\{1,2\}}$.}[.5\textwidth]{%
\begin{tikzpicture}[hide labels, scale=0.8]
\tikzstyle{VertexStyle}=[shape = circle, fill = blue!20, minimum size = 0pt, scale=0., text = white, hide labels]
\draw[black, thick] (0,0) rectangle (5,6);
\node[anchor=east,font=\small] at (0.8,5) {1};
\node[anchor=east,font=\small] at (0.8,4) {2};
\node[anchor=east,font=\small] at (0.8,3) {3};
\node[anchor=east,font=\small] at (0.8,2) {4};
\node[anchor=east,font=\small] at (0.8,1) {5};
\node[anchor=south,font=\small] at (1,0) {1};
\node[anchor=south,font=\small] at (2,0) {2};
\node[anchor=south,font=\small] at (3,0) {3};
\node[anchor=south,font=\small] at (4,0) {4};
\filldraw [gray] (1,1) circle (2pt);
\filldraw [gray] (2,1) circle (2pt);
\filldraw [gray] (3,1) circle (2pt);
\filldraw [gray] (4,1) circle (2pt);
\filldraw [gray] (1,2) circle (2pt);
\filldraw [gray] (2,2) circle (2pt);
\filldraw [gray] (3,2) circle (2pt);
\filldraw [gray] (1,3) circle (2pt);
\filldraw [gray] (2,3) circle (2pt);
\filldraw [gray] (3,3) circle (2pt);
\filldraw [gray] (4,3) circle (2pt);
\filldraw [gray] (2,3) circle (2pt);
\filldraw [gray] (1,4) circle (2pt);
\filldraw [gray] (2,4) circle (2pt);
\filldraw [gray] (1,5) circle (2pt);
\filldraw [gray] (2,5) circle (2pt);
\filldraw [gray] (3,5) circle (2pt);
\draw[very thick] (1,5) -- (1,4);
\draw[very thick] (2,5) -- (2,4);
\draw[very thick] (2,5) -- (3,5);
\draw[very thick] (1,2) -- (1,1);
\draw[very thick] (2,3) -- (2,2);
\draw[very thick] (2,1) -- (3,1) -- (4,1);
\draw[very thick] (3,2) -- (4,3);
\node (1) [label=above:{}] at (1,5) {};
\node (2) [label=above:{}] at (2,5) {};
\node (3) [label=above:{}] at (3,5) {};
\node (4) [label=above:{}] at (4,5) {};
\node (5) [label=above:{}] at (1,4) {};
\node (6) [label=above:{}] at (2,4) {};
\node (7) [label=above:{}] at (3,4) {};
\node (8) [label=above:{}] at (4,4) {};
\draw[very thick,blue] \convexpath{1,3,7,5}{.2cm};
\draw[blue,line width=1mm, ->] node[font=\fontsize{12}{0}\selectfont, right=of 2, right=1.5cm, below=-0.9cm] {$~~~~~~\rho_{\{1,2\}}$};
 \begin{pgfonlayer}{background}
    \filldraw [line width=4mm,black!3]
      (0.27,0.27)  rectangle (4.76,5.76);
  \end{pgfonlayer}
\end{tikzpicture}}%
\subcaptionbox{Splitting $\rho$ into connected subpartitions $\rho_{b_1}$, $\rho_{b_2}$.}[.5\textwidth]{%
\begin{tikzpicture}[scale=0.8] 
\draw[black, thick] (0,0) rectangle (5,6);
\node[anchor=east,font=\small] at (0.8,5) {1};
\node[anchor=east,font=\small] at (0.8,4) {2};
\node[anchor=east,font=\small] at (0.8,3) {3};
\node[anchor=east,font=\small] at (0.8,2) {4};
\node[anchor=east,font=\small] at (0.8,1) {5};
\node[anchor=south,font=\small] at (1,0) {1};
\node[anchor=south,font=\small] at (2,0) {2};
\node[anchor=south,font=\small] at (3,0) {3};
\node[anchor=south,font=\small] at (4,0) {4};
\filldraw [gray] (1,1) circle (2pt);
\filldraw [gray] (2,1) circle (2pt);
\filldraw [gray] (3,1) circle (2pt);
\filldraw [gray] (4,1) circle (2pt);
\filldraw [gray] (1,2) circle (2pt);
\filldraw [gray] (2,2) circle (2pt);
\filldraw [gray] (3,2) circle (2pt);
\filldraw [gray] (1,3) circle (2pt);
\filldraw [gray] (2,3) circle (2pt);
\filldraw [gray] (3,3) circle (2pt);
\filldraw [gray] (4,3) circle (2pt);
\filldraw [gray] (2,3) circle (2pt);
\filldraw [gray] (1,4) circle (2pt);
\filldraw [gray] (2,4) circle (2pt);
\filldraw [gray] (1,5) circle (2pt);
\filldraw [gray] (2,5) circle (2pt);
\filldraw [gray] (3,5) circle (2pt);
\draw[very thick] (1,5) -- (1,4);
\draw[very thick] (2,5) -- (2,4);
\draw[very thick] (2,5) -- (3,5);
\draw[very thick] (1,2) -- (1,1);
\draw[very thick] (2,3) -- (2,2);
\draw[very thick] (2,1) -- (3,1) -- (4,1);
\draw[very thick] (3,2) -- (4,3);
\node (1) [label=above:{}] at (1,5) {};
\node (2) [label=above:{}] at (2,5) {};
\node (3) [label=above:{}] at (3,5) {};
\node (4) [label=above:{}] at (4,3) {};
\node (5) [label=above:{}] at (1,4) {};
\node (6) [label=above:{}] at (1,3) {};
\node (7) [label=above:{}] at (3,4) {};
\node (8) [label=above:{}] at (1,1) {};
\node (9) [label=above:{}] at (4,1) {};
\node (10) [label=above:{}] at (4,2) {};
\node (11) [label=above:{}] at (4,5) {};
\draw[very thick,blue] \convexpath{1,3,7,5}{.2cm};
\draw[blue,line width=1mm, ->] node[font=\fontsize{12}{0}\selectfont, right=of 11, right=1.5cm, below=0 cm] {$~~~~~~\rho_{b_1}~$};
\draw[very thick,blue] \convexpath{6,4,9,8}{.2cm};
\draw[blue,line width=1mm, ->] node[font=\fontsize{12}{0}\selectfont, right=of 10, left=-.5cm, below=0 cm] {$~~~~~~~\rho_{b_2}~~$};
 \begin{pgfonlayer}{background}
    \filldraw [line width=4mm,black!3]
      (0.27,0.27)  rectangle (4.76,5.76);
  \end{pgfonlayer}
\end{tikzpicture}}%
\caption{Diagram $\Gamma(\rho,\pi)$ and splitting of the partition $\rho$ with $\rho\vee\pi=\{\pi_1\cup\pi_2,\pi_3\cup\pi_4\cup\pi_5\}$.}
\label{fig:diagram3}
\end{figure}

\vspace{-0.3cm}

\noindent 
 We note that for $b\subset [n]$ we have $\pi_b=\{\pi_i:i\in b\}$, 
 and any partition $\rho \in \Pi ( \pi_1\cup \cdots \cup \pi_n )$ 
 can be split into subpartitions deduced
 from the connected components of $\Gamma(\rho,\pi)$, i.e. 
 \begin{equation}
\nonumber 
   \rho=\bigcup_{b\in\sigma_\rho}\rho_b, 
\end{equation} 
 as illustrated in Figure~\ref{fig:diagram3}-$b)$ with
 $b_1 = \{ 1,2\}$, $b_2 = \{3,4,5\}$,
 and $\sigma_\rho = \{ b_1,b_2\}$. 
\begin{definition}
 For $\sigma\in\Pi([n])$ 
 we let $\Pi_{\sigma}(\pi_1\cup \cdots \cup \pi_n )$ 
 denote the collection of partitions
$\rho\in\Pi(\pi_1\cup \cdots \cup \pi_n)$ such that
$$\rho\vee\pi=\bigg\{\bigcup_{i\in b}\pi_i:b\in\sigma\bigg\}.
$$
\end{definition} 
In particular,
$\Pi_{\widehat{1}}(\pi_1\cup \cdots \cup \pi_n )$
represents the
set of connected partitions of 
$\pi_1\cup \cdots \cup \pi_n $,
and
$\Pi_{\widehat{0}}(\pi_1\cup \cdots \cup \pi_n )$
represents the partitions of 
$\pi_1\cup \cdots \cup \pi_n$
that are finer than
$\pi := \{ \pi_1,\ldots , \pi_n \}$.

\medskip
 
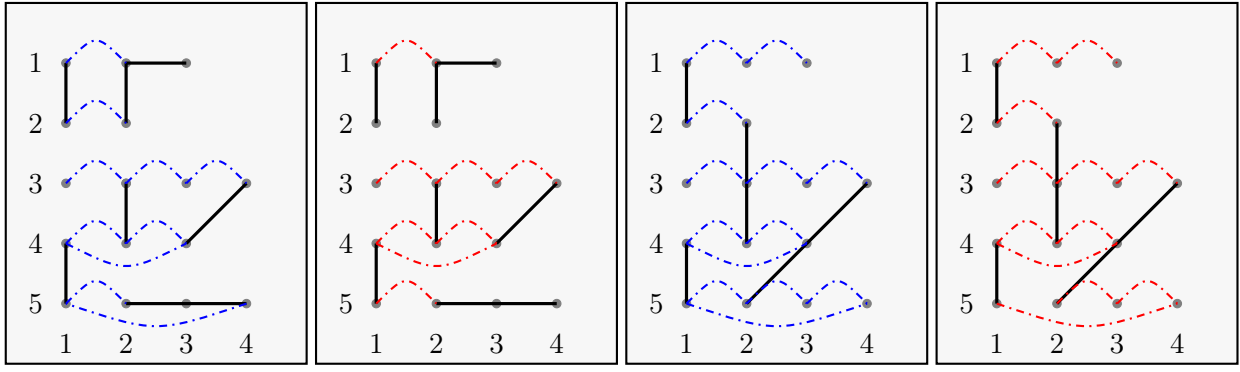
\begin{figure}[H]
\captionsetup[subfigure]{font=footnotesize}
\centering
\subcaptionbox{Multigraph $\widetilde{\rho}_G$ in blue.}[.25\textwidth]{%
\begin{tikzpicture}[hide labels,scale=0.8]
\tikzstyle{VertexStyle}=[shape = circle, fill = blue!20, minimum size = 0pt, scale=0., text = white, hide labels]
\draw[black, thick] (0,0) rectangle (5,6);
\node[anchor=east,font=\small] at (0.8,5) {1};
\node[anchor=east,font=\small] at (0.8,4) {2};
\node[anchor=east,font=\small] at (0.8,3) {3};
\node[anchor=east,font=\small] at (0.8,2) {4};
\node[anchor=east,font=\small] at (0.8,1) {5};
\node[anchor=south,font=\small] at (1,0) {1};
\node[anchor=south,font=\small] at (2,0) {2};
\node[anchor=south,font=\small] at (3,0) {3};
\node[anchor=south,font=\small] at (4,0) {4};
\filldraw [gray] (1,1) circle (2pt);
\filldraw [gray] (2,1) circle (2pt);
\filldraw [gray] (3,1) circle (2pt);
\filldraw [gray] (4,1) circle (2pt);
\filldraw [gray] (1,2) circle (2pt);
\filldraw [gray] (2,2) circle (2pt);
\filldraw [gray] (3,2) circle (2pt);
\filldraw [gray] (1,3) circle (2pt);
\filldraw [gray] (2,3) circle (2pt);
\filldraw [gray] (3,3) circle (2pt);
\filldraw [gray] (4,3) circle (2pt);
\filldraw [gray] (2,3) circle (2pt);
\filldraw [gray] (1,4) circle (2pt);
\filldraw [gray] (2,4) circle (2pt);
\filldraw [gray] (1,5) circle (2pt);
\filldraw [gray] (2,5) circle (2pt);
\filldraw [gray] (3,5) circle (2pt);
\draw[very thick] (1,5) -- (1,4);
\draw[very thick] (2,5) -- (2,4);
\draw[very thick] (2,5) -- (3,5);
\draw[very thick] (1,2) -- (1,1);
\draw[very thick] (2,3) -- (2,2);
\draw[very thick] (2,1) -- (3,1) -- (4,1);
\draw[very thick] (3,2) -- (4,3);
\node (1) [label=above:{}] at (1,5) {};
\node (2) [label=above:{}] at (2,5) {};
\node (3) [label=above:{}] at (3,5) {};
\node (4) [label=above:{}] at (4,5) {};
\node (5) [label=above:{}] at (1,4) {};
\node (6) [label=above:{}] at (2,4) {};
\node (7) [label=above:{}] at (3,4) {};
\node (8) [label=above:{}] at (4,4) {};
\draw[thick,dash dot,blue] (1,5) .. controls (1.5,5.5) .. (2,5);
\draw[thick,dash dot,blue] (1,4) .. controls (1.5,4.5) .. (2,4);
\draw[thick,dash dot,blue] (1,3) .. controls (1.5,3.5) .. (2,3);
\draw[thick,dash dot,blue] (2,3) .. controls (2.5,3.5) .. (3,3);
\draw[thick,dash dot,blue] (3,3) .. controls (3.5,3.5) .. (4,3);
\draw[thick,dash dot,blue] (1,2) .. controls (1.5,2.5) .. (2,2);
\draw[thick,dash dot,blue] (2,2) .. controls (2.5,2.5) .. (3,2);
\draw[thick,dash dot,blue] (1,2) .. controls (2,1.5) .. (3,2);
\draw[thick,dash dot,blue] (1,1) .. controls (1.5,1.5) .. (2,1);
\draw[thick,dash dot,blue] (4,1) .. controls (2.5,0.5) .. (1,1);

 \begin{pgfonlayer}{background}
    \filldraw [line width=4mm,black!3]
      (0.27,0.27)  rectangle (4.76,5.76);
  \end{pgfonlayer}
\end{tikzpicture}}%
\subcaptionbox{Graph $\rho_G$ in red.}[.25\textwidth]{
\begin{tikzpicture}[scale=0.8] 
\draw[black, thick] (0,0) rectangle (5,6);
\node[anchor=east,font=\small] at (0.8,5) {1};
\node[anchor=east,font=\small] at (0.8,4) {2};
\node[anchor=east,font=\small] at (0.8,3) {3};
\node[anchor=east,font=\small] at (0.8,2) {4};
\node[anchor=east,font=\small] at (0.8,1) {5};
\node[anchor=south,font=\small] at (1,0) {1};
\node[anchor=south,font=\small] at (2,0) {2};
\node[anchor=south,font=\small] at (3,0) {3};
\node[anchor=south,font=\small] at (4,0) {4};
\filldraw [gray] (1,1) circle (2pt);
\filldraw [gray] (2,1) circle (2pt);
\filldraw [gray] (3,1) circle (2pt);
\filldraw [gray] (4,1) circle (2pt);
\filldraw [gray] (1,2) circle (2pt);
\filldraw [gray] (2,2) circle (2pt);
\filldraw [gray] (3,2) circle (2pt);
\filldraw [gray] (1,3) circle (2pt);
\filldraw [gray] (2,3) circle (2pt);
\filldraw [gray] (3,3) circle (2pt);
\filldraw [gray] (4,3) circle (2pt);
\filldraw [gray] (2,3) circle (2pt);
\filldraw [gray] (1,4) circle (2pt);
\filldraw [gray] (2,4) circle (2pt);
\filldraw [gray] (1,5) circle (2pt);
\filldraw [gray] (2,5) circle (2pt);
\filldraw [gray] (3,5) circle (2pt);
\draw[very thick] (1,5) -- (1,4);
\draw[very thick] (2,5) -- (2,4);
\draw[very thick] (2,5) -- (3,5);
\draw[very thick] (1,2) -- (1,1);
\draw[very thick] (2,3) -- (2,2);
\draw[very thick] (2,1) -- (3,1) -- (4,1);
\draw[very thick] (3,2) -- (4,3);
\node (1) [label=above:{}] at (1,5) {};
\node (2) [label=above:{}] at (2,5) {};
\node (3) [label=above:{}] at (3,5) {};
\node (4) [label=above:{}] at (4,5) {};
\node (5) [label=above:{}] at (1,4) {};
\node (6) [label=above:{}] at (2,4) {};
\node (7) [label=above:{}] at (3,4) {};
\node (8) [label=above:{}] at (4,4) {};
\draw[thick,dash dot,red] (1,5) .. controls (1.5,5.5) .. (2,5);
\draw[thick,dash dot,red] (1,3) .. controls (1.5,3.5) .. (2,3);
\draw[thick,dash dot,red] (2,3) .. controls (2.5,3.5) .. (3,3);
\draw[thick,dash dot,red] (3,3) .. controls (3.5,3.5) .. (4,3);
\draw[thick,dash dot,red] (1,2) .. controls (1.5,2.5) .. (2,2);
\draw[thick,dash dot,red] (2,2) .. controls (2.5,2.5) .. (3,2);
\draw[thick,dash dot,red] (1,2) .. controls (2,1.5) .. (3,2);
\draw[thick,dash dot,red] (1,1) .. controls (1.5,1.5) .. (2,1);

 \begin{pgfonlayer}{background}
    \filldraw [line width=4mm,black!3]
      (0.27,0.27)  rectangle (4.76,5.76);
  \end{pgfonlayer}
\end{tikzpicture}}%
\subcaptionbox{Multigraph $\widetilde{\rho}_G$ in blue.}[.25\textwidth]{
\begin{tikzpicture}[scale=0.8] 
\draw[black, thick] (0,0) rectangle (5,6);

\node[anchor=east,font=\small] at (0.8,5) {1};
\node[anchor=east,font=\small] at (0.8,4) {2};
\node[anchor=east,font=\small] at (0.8,3) {3};
\node[anchor=east,font=\small] at (0.8,2) {4};
\node[anchor=east,font=\small] at (0.8,1) {5};

\node[anchor=south,font=\small] at (1,0) {1};
\node[anchor=south,font=\small] at (2,0) {2};
\node[anchor=south,font=\small] at (3,0) {3};
\node[anchor=south,font=\small] at (4,0) {4};

\filldraw [gray] (1,1) circle (2pt);
\filldraw [gray] (2,1) circle (2pt);
\filldraw [gray] (3,1) circle (2pt);
\filldraw [gray] (4,1) circle (2pt);
\filldraw [gray] (1,2) circle (2pt);
\filldraw [gray] (2,2) circle (2pt);
\filldraw [gray] (3,2) circle (2pt);
\filldraw [gray] (1,3) circle (2pt);
\filldraw [gray] (2,3) circle (2pt);
\filldraw [gray] (3,3) circle (2pt);
\filldraw [gray] (4,3) circle (2pt);
\filldraw [gray] (2,3) circle (2pt);
\filldraw [gray] (1,4) circle (2pt);
\filldraw [gray] (2,4) circle (2pt);
\filldraw [gray] (1,5) circle (2pt);
\filldraw [gray] (2,5) circle (2pt);
\filldraw [gray] (3,5) circle (2pt);

\draw[very thick] (1,5) -- (1,4); 

\draw[very thick] (1,2) -- (1,1);
\draw[very thick] (2,2) -- (2,4);
\draw[very thick] (2,1) -- (3,2) -- (4,3);
\draw[thick,dash dot,blue] (1,5) .. controls (1.5,5.5) .. (2,5);
\draw[thick,dash dot,blue] (2,5) .. controls (2.5,5.5) .. (3,5);
\draw[thick,dash dot,blue] (1,4) .. controls (1.5,4.5) .. (2,4);
\draw[thick,dash dot,blue] (1,3) .. controls (1.5,3.5) .. (2,3);
\draw[thick,dash dot,blue] (2,3) .. controls (2.5,3.5) .. (3,3);
\draw[thick,dash dot,blue] (3,3) .. controls (3.5,3.5) .. (4,3);
\draw[thick,dash dot,blue] (1,2) .. controls (1.5,2.5) .. (2,2);
\draw[thick,dash dot,blue] (2,2) .. controls (2.5,2.5) .. (3,2);
\draw[thick,dash dot,blue] (1,2) .. controls (2,1.5) .. (3,2);
\draw[thick,dash dot,blue] (1,1) .. controls (1.5,1.5) .. (2,1);
\draw[thick,dash dot,blue] (2,1) .. controls (2.5,1.5) .. (3,1);
\draw[thick,dash dot,blue] (3,1) .. controls (3.5,1.5) .. (4,1);
\draw[thick,dash dot,blue] (1,1) .. controls (2.5,0.5) .. (4,1);

 \begin{pgfonlayer}{background}
    \filldraw [line width=4mm,black!3]
      (0.27,0.27)  rectangle (4.76,5.76);
  \end{pgfonlayer}
\end{tikzpicture}}%
\subcaptionbox{Graph and $\rho_G$ in red.}[.25\textwidth]{
\begin{tikzpicture}[scale=0.8] 
\draw[black, thick] (0,0) rectangle (5,6);

\node[anchor=east,font=\small] at (0.8,5) {1};
\node[anchor=east,font=\small] at (0.8,4) {2};
\node[anchor=east,font=\small] at (0.8,3) {3};
\node[anchor=east,font=\small] at (0.8,2) {4};
\node[anchor=east,font=\small] at (0.8,1) {5};

\node[anchor=south,font=\small] at (1,0) {1};
\node[anchor=south,font=\small] at (2,0) {2};
\node[anchor=south,font=\small] at (3,0) {3};
\node[anchor=south,font=\small] at (4,0) {4};

\filldraw [gray] (1,1) circle (2pt);
\filldraw [gray] (2,1) circle (2pt);
\filldraw [gray] (3,1) circle (2pt);
\filldraw [gray] (4,1) circle (2pt);
\filldraw [gray] (1,2) circle (2pt);
\filldraw [gray] (2,2) circle (2pt);
\filldraw [gray] (3,2) circle (2pt);
\filldraw [gray] (1,3) circle (2pt);
\filldraw [gray] (2,3) circle (2pt);
\filldraw [gray] (3,3) circle (2pt);
\filldraw [gray] (4,3) circle (2pt);
\filldraw [gray] (2,3) circle (2pt);
\filldraw [gray] (1,4) circle (2pt);
\filldraw [gray] (2,4) circle (2pt);
\filldraw [gray] (1,5) circle (2pt);
\filldraw [gray] (2,5) circle (2pt);
\filldraw [gray] (3,5) circle (2pt);

\draw[very thick] (1,5) -- (1,4); 

\draw[very thick] (1,2) -- (1,1);
\draw[very thick] (2,2) -- (2,4);
\draw[very thick] (2,1) -- (3,2) -- (4,3);
\draw[thick,dash dot,red] (1,5) .. controls (1.5,5.5) .. (2,5);
\draw[thick,dash dot,red] (2,5) .. controls (2.5,5.5) .. (3,5);
\draw[thick,dash dot,red] (1,4) .. controls (1.5,4.5) .. (2,4);
\draw[thick,dash dot,red] (1,3) .. controls (1.5,3.5) .. (2,3);
\draw[thick,dash dot,red] (2,3) .. controls (2.5,3.5) .. (3,3);
\draw[thick,dash dot,red] (3,3) .. controls (3.5,3.5) .. (4,3);
\draw[thick,dash dot,red] (1,2) .. controls (1.5,2.5) .. (2,2);
\draw[thick,dash dot,red] (2,2) .. controls (2.5,2.5) .. (3,2);
\draw[thick,dash dot,red] (1,2) .. controls (2,1.5) .. (3,2);
\draw[thick,dash dot,red] (2,1) .. controls (2.5,1.5) .. (3,1);
\draw[thick,dash dot,red] (3,1) .. controls (3.5,1.5) .. (4,1);
\draw[thick,dash dot,red] (1,1) .. controls (2.5,0.5) .. (4,1);

 \begin{pgfonlayer}{background}
    \filldraw [line width=4mm,black!3]
      (0.27,0.27)  rectangle (4.76,5.76);
  \end{pgfonlayer}
\end{tikzpicture}}%
\caption{Diagram $\Gamma ( \rho , \pi )$,
  multigraph $\widetilde{\rho}_G$, and graph $\rho_G$.}
\label{fig:diagram5}
\end{figure}

 Given $F:\Pi'(\pi_1\cup \cdots \cup \pi_n )\to\R$, 
 where $\Pi'(\pi_1\cup \cdots \cup \pi_n )$ is the collection of all subpartitions of $\pi_1\cup \cdots \cup \pi_n $, 
 we define the mixed moments $\widehat{F}:2^{[n]}\to\R$ by 
\begin{equation}
\label{mm} 
  \widehat{F}(A)=\sum_{\rho\in\Pi(\cup_{i\in A}\pi_i)}F(\rho),
  \qquad
   A\subset [n], 
\end{equation}
cf. \cite[p.~33]{MalyshevMinlos91}.
The semi-invariants $C_F:2^{[n]}\to\R$ are defined by the induction formula
 $C_F(A)=\widehat{F}(A)$ when $|A|=1$, and 
\begin{equation}
\nonumber
C_F(A)=\widehat{F}(A)-\sum_{\substack{\{b_1,\dots,b_k\}\in\Pi(A)\\k\ge2}}\prod_{i=1}^kC_F(b_i),
\end{equation}
for $|A|>1$,
see Relation~(16) page~33 of \cite{MalyshevMinlos91}, 
where the sum is taken over all partitions $\sigma\in\Pi(A)$
such that $|\sigma|\ge2$,
i.e.
\begin{equation}
\label{2} 
C_F(A)=\sum_{\rho\in\Pi(A)}
  (-1)^{|\rho|} (|\rho|-1)!
 \prod_{b\in\rho} \widehat{F}( b),
\end{equation}
 see Relation~(16') in \cite{MalyshevMinlos91}.
 The next proposition generalizes \cite[Proposition~3.3]{LiuPrivault} to
the multivariate case.
\begin{prop}
 Suppose that $F$ satisfies the connectedness factorization property 
\begin{equation}\label{factor-1}
  F(\rho)=\prod_{b\in\sigma_\rho}F(\rho_b),
  \qquad
   \rho\in\Pi'(\pi_1\cup \cdots \cup \pi_n ). 
\end{equation}
 Then, the semi-invariants are given by
 \begin{equation}
   \label{eq-1} 
  C_F(A)=\sum_{\rho\in\Pi_{\widehat{1}}(\cup_{i\in A}\pi_i)}F(\rho),
  \qquad
  \emptyset \not= A\subset [n].
\end{equation}
\end{prop}
\begin{Proof}
 \noindent
 $(i)$ 
 It is clear that \eqref{eq-1} holds when $|A|=1$. 
 When $|A|=2$, taking $A=\{i,j\} \subset [n]$,
 $i\not= j$, we have 
\begin{eqnarray*}
C_F(A)&=&\widehat{F}(\{i,j\})-C_F(\{i\})C_F(\{j\})\\
&=&\sum_{\rho\in\Pi(\pi_i \cup\pi_j )}F(\rho)-\widehat{F}(\{i\})\widehat{F}(\{j\})\\
&=&\sum_{\rho\in\Pi_{\widehat{1}}(\pi_i \cup\pi_j )}F(\rho)+\sum_{\rho\in\Pi_{\widehat{0}}(\pi_i \cup\pi_j )}F(\rho)-\left(
\sum_{\rho_1\in\Pi(\pi_i)}F(\rho_1)\right)
\left(
\sum_{\rho_2\in\Pi(\pi_j)}F(\rho_2)\right).
\end{eqnarray*}
By splitting any $\rho\in\Pi_{\widehat{0}}(\pi_i\cup\pi_j)$
into two disjoint subpartitions according to Definition~\ref{part-1}, i.e. 
$$\rho=\rho_{\{i\}}\cup\rho_{\{j\}}, 
$$
together with the factorization property \eqref{factor-1}, we find 
\begin{eqnarray}
  \nonumber
  \sum_{\rho\in\Pi_{\widehat{0}}(\pi_i \cup\pi_j)}F(\rho)&=&\sum_{\substack{\rho\in\Pi_{\widehat{0}}(\pi_i\cup\pi_j)
      \\
\nonumber
\rho=\rho_{\{i\}}\cup\rho_{\{j\}}}}F(\rho_{\{i\}})F(\rho_{\{j\}})
  \\
\nonumber
  &=&\left(
\sum_{\rho_1\in\Pi(\pi_i)}F(\rho_1)\right)
\left(
\sum_{\rho_2\in\Pi(\pi_j)}F(\rho_2)\right), 
\end{eqnarray}
 which shows \eqref{eq-1}. 

 \smallskip

 \noindent
 $(ii)$ 
 Next, suppose that \eqref{eq-1} holds for any $A\subset [n]$ with $|A|\leq l \leq n-1$. Let $A\subset [n]$ be a subset of $[n]$ with $|A|=l +1$. We have 
\begin{eqnarray*}
  \widehat{F}(A)&=&
    \sum_{\substack{\rho\in\Pi(\cup_{i\in A}\pi_i)
          }
    }F(\rho)
  \\
  &=&
  \sum_{\substack{\sigma=\{b_1,\dots,b_k\}\in\Pi(A)\\k\geq 1}}\sum_{\substack{\rho\in\Pi(\cup_{i\in A}\pi_i)\\
\rho\vee\pi_A=\{\cup_{i\in b_j}\pi_i\}_{j=1}^k}}F(\rho)\\
  &=&
  \sum_{\substack{\sigma=\{b_1,\dots,b_k\}\in\Pi(A)\\k\geq 1}}\sum_{\substack{\rho\in\Pi(\cup_{i\in A}\pi_i)\\
\rho\vee\pi_A=\{\cup_{i\in b_j}\pi_i\}_{j=1}^k}}\prod_{j=1}^kF(\rho_{b_j})\\
  &=&
  \sum_{\substack{\sigma=\{b_1,\dots,b_k\}\in\Pi(A)\\k\geq 1}}\prod_{j=1}^k
 \sum_{\substack{\rho_j\in\Pi(\cup_{i\in b_j}\pi_i)\\
     \rho_j\vee\pi_{b_j}=\widehat{1}}}F(\rho_{b_i})
 \\
 &=&
 \sum_{\substack{\sigma=\{b_1,\dots,b_k\}\in\Pi(A)\\k\geq 1}}\prod_{j=1}^k
\sum_{\rho_j\in\Pi_{\widehat{1}}(\cup_{i\in b_j}\pi_i)}F(\rho_{b_i})
\\
&=& \sum_{\rho\in\Pi_{\widehat{1}}(\cup_{i\in A}\pi_i)}F(\rho)
+ \sum_{\substack{\{b_1,\dots,b_k\}\in\Pi(A)\\k\geq 2}}\prod_{j=1}^kC_F(b_j),
\end{eqnarray*}
where the last equality follows from the induction hypothesis \eqref{eq-1} when $|A|\leq l $. The proof is completed by subtracting the last term from both sides.
\end{Proof}
 Given $n\ge1$ and $f^{(i)}:(\R^d)^{r_i}\to\R$, $i=1,\dots,n$,
 measurable functions, we let 
$$
 \left(\bigotimes_{i=1}^nf^{(i)} \right)(x_{1,1},\ldots,x_{1,r_1},\ldots,
 x_{n,1},\ldots , x_{n,r_n}):=\prod_{i=1}^nf^{(i)} (x_{i,1},\dots,x_{i,r_i}). 
  $$
  For $\rho\in \Pi(\pi_1\cup \cdots \cup \pi_n )$, we also denote by
  $\big(
  \bigotimes_{i=1}^nf^{(i)} \big)_{\hskip-0.01cm \rho}:(\R^d)^{|\rho|}\to\R$ the function
  obtained by equating any two variables
  whose indexes belong to a same block of $\rho$.
  We refer to \cite[Theorem~3.1]{bogdan}
  for the next result.
    \begin{prop}
    \label{moment-1}
    Let $n\ge1$,
    $r_1,\ldots ,r_n\ge1$,
    and let $f^{(i)} :(\R^d)^{r_i}\to\R$ be a
    sufficiently integrable measurable function 
    for $i=1,\dots,n$. We have
\begin{equation}
\nonumber
\E\left[\prod_{i=1}^n
 \sum_{(x_1, \ldots ,x_{r_i})\in\eta^{r_i} }
 f^{(i)} (x_1,\dots,x_{r_i})
 \right]=\sum_{\rho\in\Pi(\pi_1\cup \cdots \cup \pi_n )}\lambda^{|\rho|}\int_{(\R^d)^{|\rho|}}\left(\bigotimes_{i=1}^n f^{(i)} \right)_{\hskip-0.1cm \rho} \! (\mathbf{x})
 \ \! \mu^{\otimes |\rho|} ( \mathrm{d}\mathbf{x} ),
\end{equation}
\end{prop}
Proposition~\ref{moment-1} can be specialized as follows.
\begin{corollary}\label{moment-2}
  Let $r_i\ge2$, $i=1,\dots,n$, 
  and consider $f^{(i)}:(\R^d)^{r_i}\to\R$
  measurable functions that vanish on diagonals,
  i.e. $f^{(i)} (x_1,\dots,x_{r_i})=0$ whenever $x_k=x_l$ for some $1\leq k\neq l\leq r_i$, $i=1,\ldots , n$. We have   
\begin{equation}
  \label{fjkldf4}
  \E\left[\prod_{i=1}^n
    \sum_{(x_1, \ldots ,x_{r_i})\in\eta^{r_i} }
    f^{(i)} (x_1,\dots,x_{r_i})
\right]= 
  \sum_{\substack{\rho\in\Pi ( \pi_1 \cup \cdots \cup \pi_n ) 
        \\\rho\wedge\pi=\widehat{0}} \atop {\rm (non-flat) \atop
  }}
  \lambda^{|\rho|}\int_{(\R^d)^{|\rho|}}\left(\bigotimes_{i=1}^n f^{(i)} \right)_{\hskip-0.1cm \rho} \! (\mathbf{x})
  \ \! \mu^{\otimes |\rho|} ( \mathrm{d}\mathbf{x} ).
\end{equation}
\end{corollary}
For $i=1,\dots,n$, let $M_i \subset \{1,\ldots , m\}$,
 $r_i\geq 2$,
 and let $G_i=(V_{G_i},E_{G_i})$
 be a connected graph  with edge set $E_{G_i}$ and
 vertex set of the form 
 $V_{G_i}=\big(v^{(i)}_1, \ldots ,v^{(i)}_{r_i}; {\{w_j^{(i)}\}_{j\in M_i}}\big)$,
 such that 
 \begin{enumerate}[i)]
 \item the subgraph $\GG_i$ induced by $G_i$ on
   $\{v^{(i)}_1, \ldots ,v^{(i)}_{r_i}\}$ is connected, and 
\item 
  {the endpoint vertices $
    \{ w_j^{(i)}\}_{j\in M_i}$} are not adjacent to each other in $G_i$, 
\end{enumerate}
 and let $G:=\{G_1,\ldots , G_n\}$.
 In Definition~\ref{fjmklc3},
 for every $\rho\in\Pi(\pi_1\cup \cdots \cup \pi_n )$
 we build a graph structure induced by
 $(G_1,\dots,G_n)$ on the diagram $\Gamma(\rho,\pi)$,
 analogous to \cite[Definition~2.2]{LiuPrivault}.
\begin{definition}
   \label{fjmklc3} 
 Given $\rho \in\Pi(\pi_1\cup \cdots \cup \pi_n )$
 a partition of $\pi_1\cup \cdots \cup \pi_n $, 
 we let $\widetilde{\rho}_G$ denote the multigraph 
 constructed as follows on $[m] \cup \pi_1 \cup \cdots \cup \pi_n$: 
\begin{enumerate}[i)]  
\item for all $j_1, j_2\in [r_i]$, $j_1\not= j_2$, and $i\in [n]$, 
  an edge links $(i,j_1)$ to $(i,j_2)$
  iff $\{v^{(i)}_{j_1},v^{(i)}_{j_2}\}\in E_{G_i}$. 
\item for all $(j,k)\in [r_i]\times M_i$ and $i\in [n]$, an edge
  links $(k)$ to $(i,j)$ iff $\{v^{(i)}_j,w_k^{(i)}\}\in E_{G_i}$; 
\item for all $i_1,i_2\in [n]$ and 
  $(j_1,j_2) \in [r_{i_1}]\times [r_{i_2}]$,
  we merge any two nodes $(i_1,j_1)$ and $(i_2,j_2)$ 
  if they belong to a same block in $\rho$. 
\end{enumerate}
In addition, we let $\rho_G$ be the graph constructed
from the multigraph $\widetilde{\rho}_G$ 
by removing any redundant edge
in $\widetilde{\rho}_G$. 
\end{definition}
As in Section~\ref{diagramrepresentation},
 the graph $\rho_G$ forms a connected graph with
 $|\rho | + m$ vertices. 
 Figure~\ref{fig:diagram5} presents two examples of
 multigraphs $\widetilde{\rho}_G$ and graphs $\rho_G$ when $G_1,G_2,G_3$ are line graphs, $G_4$ is a triangle, and $G_5$ is a rectangle 
 on a partition diagram $\Gamma ( \rho , \pi )$ with
 no endpoints, i.e. $M_1=\cdots = M_n = \emptyset$ here.
 
\vspace{-0.4cm}
\noindent
For each $i=1,\ldots , n$
 denote by $N_{M_i}^{G_i}$ the count of subgraphs in
 the random-connection model $G_H (\eta \cup \{y_j\}_{j\in M_i} )$
 with endpoint set $\{y_j\}_{j\in M_i}$, i.e.
\begin{equation}
\nonumber
N_{M_i}^{G_i}
=\sum_{(x_1, \ldots ,x_{r_i})\in\eta^{r_i}} f^{(i)}_{M_i} (x_1, \ldots ,x_{r_i}), 
\end{equation}
 where
 $f^{(i)}_{M_i} :(\real^d)^{r_i} \to \{0,1\}$ is the random function defined as 
\begin{equation}
\nonumber
f^{(i)}_{M_i} (x_1, \ldots ,x_{r_i}):=
\prod_{  \substack{
    1 \leq l \leq r_i, \ j \in M_i
    \\ \{w_j^{(i)} , v^{(i)}_l\}\in E_{G_i} }
}
\bone_{\{y_j\leftrightarrow x_l \}} 
\prod_{\substack{ 1 \leq k,l \leq r_i
    \\ \{v^{(i)}_k,v^{(i)}_l\}\in E_{G_i}}}\bone_{\{x_k \leftrightarrow x_l \}},
\qquad
 x_1,\ldots , x_{r_i} \in \R^d. 
\end{equation} 
For $\rho = \{ b_1,\ldots , b_{|\rho |}\} \in
\Pi (\pi_1 \cup \cdots \cup \pi_n )$,
we also let 
\begin{equation}
\nonumber
    {\cal A}^\rho_j:=\{ i \in [ |\rho | ] \ : \ \exists (s,k)\in b_i ~\mathrm{s.t.}~
    \big(v^{(s)}_k,w_j^{(s)}\big) \in E_{G_s} 
    \} 
\end{equation} 
denote the neighborhood of the vertex $(|\rho | + j)$ in $\rho_G$,
$j=1,\ldots , m$.
{The next proposition is a consequence of 
 Relation~\eqref{eq-1} and Corollary~\ref{moment-2}. 
} 
\begin{prop}
  \label{fjklf2}
  Let $N_{M_i}^{G_i}$ be subgraph counts in the
  random-connection model $G_H (\eta\cup\{y_1,\dots,y_m\})$
  as defined above, for $i=1,\dots,n$. We have 
    \begin{equation}\label{moment-4}
      \E\left[\prod_{i=1}^n N^{G_i}_{M_i}\right]=\sum_{
        \substack{
          \rho\in\Pi(\pi_1 \cup \cdots \cup \pi_n )
      \\\rho\wedge\pi=\widehat{0}} \atop {\rm (non-flat)}}
\lambda^{|\rho|}\int_{(\R^d)^{|\rho|}}\prod_{\substack{ 
      1 \leq j \leq m 
      \\ i\in {\cal A}^\rho_j}}
    H(x_i,y_j)\prod_{(k,l)\in E_{\rho_G} }H(x_k,x_l) \ \! \mu(\mathrm{d}\mathbf{x}),
\end{equation}
and joint cumulant 
\begin{equation}
  \label{cum-2}
  \kappa(N_{G_1},\dots,N_{G_n})=\sum_{
    \substack{\rho\in\Pi_{\widehat{1}}(\pi_1 \cup \cdots \cup \pi_n )
        \\\rho\wedge\pi=\widehat{0}} \atop {\rm (non-flat \ \! connected)}}
  \lambda^{|\rho|}\int_{(\R^d)^{|\rho|}}
  \prod_{\substack{ 
       1 \leq j \leq m 
      \\ i\in {\cal A}^\rho_j}}
    H(x_i,y_j)\prod_{(k,l)\in E_{\rho_G} }H(x_k,x_l) \ \! \mu(\mathrm{d}\mathbf{x}). 
\end{equation}
\end{prop}
\begin{Proof}
 The moment identity \eqref{moment-4} is obtained by
 taking expectation on both sides of
 \eqref{fjkldf4} in Corollary~\ref{moment-2}
 and using the relation $H(x,y) = \E [ 
 \bone_{\{x \leftrightarrow y\}}]$,
 $x, y \in \real^d$.
 Next, we note that 
 the connectedness factorization property
 \eqref{factor-1} is satisfied by
 $$F(\rho ):=\lambda^{|\rho |}\int_{(\R^d)^{|\rho |}}
 \prod_{\substack{ 
       1 \leq j \leq m 
      \\ i\in {\cal A}^\rho_j}}
 H(x_i,y_j)
 \prod_{(k,l)\in E_{\rho_G} }H(x_k,x_l) \ \! \mu(\mathrm{d}x_1)\cdots\mu(\mathrm{d}x_{|\rho |}),
 $$
 $\rho \in\Pi(\pi_1\cup \cdots \cup \pi_n )$,
 hence \eqref{cum-2} follows from
 Relations~\eqref{mm}, \eqref{eq-1}, \eqref{moment-4},  
 and the classical cumulant-moment relationship \eqref{2},
 see e.g. Relation~(3.3) in \cite{elukacs}. 
\end{Proof}
\noindent
The cumulant formula of Proposition~\ref{fjklf2}
is implemented in the code listed in Appendix~\ref{fjkldsf-2}. 
\section{Cumulant and factorial moment estimates} 
\label{statuleviciuscond}
\noindent
The following result can be found in
 \cite[Corollary~2.1]{saulis} or \cite[Theorem~2.4]{doering}. 
\begin{lemma}\label{Statuleviciuscond1}
  Let $\{X_\lambda\}$ be a family of random variables with moments
  of all orders, mean zero and unit variance for all $\lambda>0$.
  Suppose that for all $j\geq 3$ and sufficiently large $\lambda$,
  the cumulant of order $j$ of $X_\lambda$ is bounded by
\begin{equation}
\nonumber 
|\kappa_j(X_\lambda)|\leq \frac{(j!)^{1+\gamma}}{(\Delta_\lambda)^{j-2}}
\end{equation}
where $\gamma\ge0$ is a constant independent of $\lambda$.
Then we have the Berry-Esseen bound 
\begin{equation}
\nonumber
  \sup_{x\in\R}|\IP (X_\lambda\leq x)-\Phi(x)|\leq C_\gamma (\Delta_\lambda)^{-1/(1+2\gamma)},
\end{equation}
 for $C_\gamma>0$ a constant depending only on $\gamma$. 
\end{lemma}
 We let
 $m_n(X) := \E [ X(X-1) \cdots (X-n+1)]$   
 denote the factorial moments of order $n \geq 1$ 
 of a discrete random variable $X$. 
\begin{prop}[Corollary 1.13 in \cite{bollobas}] 
  \label{fdshkf0}
  Assume that 
  $$
  \lim_{n\to \infty} m_n(X)\frac{n^k}{n!} =0,
  \qquad
 k \geq 0.
  $$ 
 Then for any $n\geq 0$, we have 
 \begin{equation}
   \label{fjkl32}
   \P (   X = n ) = 
\frac{1}{n!}
\sum_{i\geq 0} \frac{(-1)^i}{i!}
m_{n + i}(X).
\end{equation}
\end{prop}
\section{Gram-Charlier expansions}
\label{s5}
\noindent
 Let $\varphi(x) := e^{-x^2/2} / \sqrt{2\pi}$,
 $x\in \real$,
 denote the standard normal probability density function. 
 In addition to the second order expansion Gaussian approximation 
\begin{equation} 
\label{gram_charlier-1}
\phi_X^{(1)}(x)=
\frac{1}{\sqrt{\kappa_2}}
\varphi \left( \frac{x-\kappa_1}{\sqrt{\kappa_2}}\right)
\end{equation}
for the probability density $\phi_X(x)$
function of a random variable $X$,
higher order Gram-Charlier expansions of 
 third and fourth order are given by 
\begin{equation} 
\label{gram_charlier-2}
\phi_X^{(3)}(x)=
\frac{1}{\sqrt{\kappa_2}}
\varphi \left( \frac{x-\kappa_1}{\sqrt{\kappa_2}}\right)
\left( 1 +
c_3 H_3\left(
\frac{x-\kappa_1}{\sqrt{\kappa_2}} \right)
\right)
\end{equation}
 and
\begin{equation} 
\nonumber 
\phi_X^{(4)}(x)=
\frac{1}{\sqrt{\kappa_2}}
\varphi \left( \frac{x-\kappa_1}{\sqrt{\kappa_2}}\right)
\left( 1 +
c_3 H_3\left(
\frac{x-\kappa_1}{\sqrt{\kappa_2}} \right)
+
c_4 H_4\left(
\frac{x-\kappa_1}{\sqrt{\kappa_2}} \right)
 + c_6 H_6\left( \frac{x-\kappa_1}{\sqrt{\kappa_2}} \right)
 \right).
\end{equation}
 see \S~17.6 of \cite{cramer}, where 
\begin{itemize}
\item 
$H_0(x)=1$, 
$H_1(x)=x$, 
$H_3(x)=x^3-3x$,
$H_4(x)=x^4-6x^2+3$, 
$H_6(x)=x^6-15x^4+45x^2-15$ 
 are Hermite polynomials, 
 \item 
  the sequence $c_3,c_4,c_5,c_6$ is given from the cumulants $(\kappa_n)_{n\geq 1}$
of $X$ as 
$$ 
c_3 = \frac{\kappa_3}{3! (\kappa_2)^{3/2}}, 
\quad
c_4 = \frac{\kappa_4}{4! (\kappa_2)^2}, 
\quad 
c_5 = \frac{\kappa_5}{5! \kappa_5^{5/2}},
\quad
 c_6 =
 \frac{\kappa_6}{6! (\kappa_2)^3}
 +
  \frac{(\kappa_3)^2}{2(3!)^2 (\kappa_2)^3}, 
$$
 where $c_3$ and $c_4$ are expressed from 
 the skewness $\kappa_3/(\kappa_2)^{3/2}$ and
 the excess kurtosis $\kappa_4/(\kappa_2)^2$.
\end{itemize} 
\noindent

\section{Cumulant code}
\label{fjkldsf}
\noindent  
The following code generates closed-form cumulant expressions 
via symbolic calculations in SageMath for any dimension $d\geq 1$,
any connected subgraph
$\GG$ induced by $G$, and any set of endpoint connections 
represented as the sequence 
$\EE = [\EE_1,\ldots , \EE_m]$.
When $G$ has no endpoint ($m=0$)
we have $\EE =[ \ ]$, 
 however, in this case the measure $\mu$ should be finite, i.e., 
the density function
\EscVerb{mu(x,λ,β)} should be
integrable with respect to the Lebesgue measure on $\real^d$.
The choice of SageMath for this implementation
is due to its fast handling of symbolic integration via Maxima, which is significantly faster than the Python package Sympy. 
This code is also sped up by parallel processing that
 can distribute the load among different CPU cores. 
This SageMath code and the next one
are available for download at 
\url{https://github.com/nprivaul/random-connection}.  

\bigskip

\begin{lstlisting}
from time import time
import datetime
import multiprocessing as mp

global cumulants

def partitions(points):
    if len(points) == 1:
        yield [ points ]
        return
    first = points[0]
    for smaller in partitions(points[1:]):
        for m, subset in enumerate(smaller):
            yield smaller[:m] + [[ first ] + subset]  + smaller[m+1:]
        yield [ [ first ] ] + smaller

def nonflat(partition,r):
    p = []
    for j in partition:    
        seq = list(map(lambda x: (x-1)//r,j))
        p.append(len(seq) == len(set(seq)))
    return all(p)

def connected(partition,n,r):
    q = []; c = 0
    if n  == 1: return all([len(j)==1 for j in partition])
    for j in partition:
        jk = list(set(map(lambda x: (x-1)//r,j)))
        if(len(jk)>1):            
            if c == 0:
                q = jk; c += 1
            elif(set(q) & set(jk)):
                d=[y for y in (q+jk) if y not in q]
                q = q + d
    return n == len(set(q))

def connectednonflat(n,r):
    points = list(range(1,n*r+1))
    randd = []
    for m, p in enumerate(partitions(points), 1):
        randd.append(sorted(p))
    cnfp = [e for e in randd if (connected(e,n,r) and nonflat(e,r))]
    for rou in range(r,(r-1)*n+2): 
        rs = [d for d in cnfp if len(d)==rou]
        print("Connected non-flat partitions with",rou,"blocks:",len(rs))
    print("Connected non-flat set partitions:",len(cnfp))
    return cnfp

def graphs(G,EP,setpartition,n):
    r=len(set(flatten(G)));rhoG = []
    for j in range(n):
        for hop in G: rhoG.append([r*j+hop[0],r*j+hop[1]])
        for l in range(len(EP)):
            F=EP[l]
            for i in F: rhoG.append([j*r+i,n*r+l+1]);
    for i in setpartition:
        if(len(i)>1):
            b = []
            for j in rhoG:
                b.append([i[0] if ele in i else ele for ele in j])
            rhoG = b
    for i in rhoG: i.sort()
    return rhoG

def inner(n,d,G,EP,mu,H,setpartition,z,r):
    rhoG=graphs(G,EP,setpartition,n)
    for ll in range(len(EP)+1):
        for l in range(1,d+1): z[d*(n*r+ll)+l] = var(str(y)+str(ll)+str('_')+str(l))
    for key in range(1,n*r+1): 
        for l in range(1,d+1): z[key*d+l] = var(str(x)+str(key)+str(x)+str(l))
    edgesrhoG = [i for n, i in enumerate(rhoG) if i not in rhoG[:n]]
    vertrhoG = set(flatten(edgesrhoG));
    for ll in range(len(EP)): vertrhoG.remove(n*r+ll+1);
    strr = '*λ'*len(vertrhoG)
    for i in vertrhoG:
        for l in range(1,d+1): strr = '*mu({},{},{})'.format(z[i*d+l], λ, β) + strr
        for l in range(1,d+1): strr = strr + ').integrate({},-infinity,+infinity)'.format(z[i*d+l])
    for i in edgesrhoG:
        for l in range(1,d+1): strr = '*H({},{},{})'.format(z[i[0]*d+l],z[i[1]*d+l], β) + strr
    strr = '('*len(vertrhoG)*d+strr[1:]
    return eval(preparse(strr))

def collect_result(result):
    global cumulants
    global iii
    global tim
    iii=iii+1;
    if (mod(iii,100)==0):
        tim=(time()-t_start2)*(lencnfp-iii)/iii/60
        print('[%d]\r'%(iii),'Est. remaining time (minutes):%d'%(tim),end="")
    cumulants+=result
	
def c(n,d,G,EP,mu,H):
    global cumulants
    global iii
    global t_start2
    t_start2 = time()
    d_start2 = datetime.datetime.now()
    r=len(set(flatten(G)));
    x,y=var("x,y")
    cumulants = 0; iii = 0
    z = dict(enumerate([str(x)+str(key)+str(x)+str(l) for key in range(0,n*r+1) for l in range(1,d+1)], start=1))
    global lencnfp
    cnfp=connectednonflat(n,r)
    lencnfp=len(cnfp)
    pool = mp.Pool(4) # pool = mp.Pool(mp.cpu_count())
    for setpartition in cnfp: 
        pool.apply_async(func = inner, args=(n,d,G,EP,mu,H,setpartition,z,r), callback=collect_result)
    pool.close()
    pool.join()
    print("\n");
    d_end2 = datetime.datetime.now()
    print("Runtime is",(d_end2-d_start2))
    return cumulants._sympy_() 
\end{lstlisting}

\section{Joint cumulant code}
\label{fjkldsf-2}
\noindent  
The following code generates closed-form joint cumulant expressions 
via symbolic calculations in SageMath for any dimension $d\geq 1$,
any sequence $(\GG_1,\ldots , \GG_n)$ of connected subgraphs
induced by $(G_1,\ldots ,G_n)$. 
As above, the endpoint connections are represented using 
represented as the sequence 
$\EE = [\EE_1,\ldots , \EE_m]$
where $\EE_i$ denotes the set of vertices of $\GG_1\cup \cdots \cup \GG_n$
which are attached to the $i$-$th$
endpoint, $i=1,\ldots , m$. 

\medskip

\begin{lstlisting}
def jpartitions(points):
    if len(points) == 1:
        yield [ points ]
        return
    first = points[0]
    for smaller in jpartitions(points[1:]):
        for m, subset in enumerate(smaller):
            yield smaller[:m] + [[ first ] + subset]  + smaller[m+1:]
        yield [ [ first ] ] + smaller

def jnonflat(partition,rr):
    n=len(rr); p = []
    for j in partition:    
        for i in range(n):
            j2 = [l for l in j if l > sum(rr[0:i]) and l<=sum(rr[0:(i+1)])]
            p.append(len(j2) <= 1)
    return all(p)

def jconnected(partition,rr):
    n=len(rr); q = []; c = 0; 
    if n  == 1: return True
    for j in partition:
        jk = [i for i in range(n) if len([l for l in j if l > sum(rr[0:i]) and l<=sum(rr[0:(i+1)])])>=1]
        if(len(jk)>1):            
            if c == 0:
                q = jk; c += 1
            elif(set(q) & set(jk)):
                d=[y for y in (q+jk) if y not in q]
                q = q + d
    return n == len(q)

def jconnectednonflat(rr):
    n=len(rr); 
    points = list(range(1,sum(rr)+1))
    randd = []
    for m, p in enumerate(jpartitions(points), 1): randd.append(sorted(p))
    for rou in range(min(rr),sum(rr)-n+2):    
        rs = [d for d in randd if (jnonflat(d,rr) and len(d)==rou)]
        rss = [e for e in rs if jconnected(e,rr)]
        print("Connected non-flat partitions with",rou,"blocks:",len(rss))
    cnfp = [e for e in randd if (jconnected(e,rr) and jnonflat(e,rr))]
    print("Connected non-flat set partitions:",len(cnfp))
    return cnfp

def jgraphs(G,EP,setpartition):
    rr=[len(set(flatten(g))) for g in G];
    n=len(G); rhoG = []
    ee=[len(set(flatten(e))) for e in EP];
    for j in range(n):
        for hop in G[j]: rhoG.append([hop[0],hop[1]])
        for l in range(len(EP)):
            F=EP[l]
            for i in F: rhoG.append([i,sum(rr)+l+1]);
    for i in setpartition:
        if(len(i)>1):
            b = []
            for j in rhoG:
                b.append([i[0] if ele in i else ele for ele in j])
            rhoG = b
    for i in rhoG: i.sort()
    return rhoG

def jc(d,G,EP,mu,H):
    rr=[len(set(flatten(g))) for g in G];
    if(sum(rr)!=len(set(flatten(G)))):
        print("Wrong G format");
        return 0
    n=len(G);
    ee=[len(set(flatten(e))) for e in EP];
    x,y=var("x,y")
    jcumulants = 0; ii=0
    z = dict(enumerate([str(x)+str(key)+str(x)+str(l) for key in range(1,sum(rr)+1) for l in range(1,d+1)], start=1))
    cnfp=jconnectednonflat(rr)
    for setpartition in cnfp: 
        ii=ii+1;print('[%d/%d]\r'%(ii,len(cnfp)),end="")
        rhoG=jgraphs(G,EP,setpartition)
        for j in range(n):
            m=len(EP); 
            for l in range(m+1):
                for ld in range(1,d+1): z[d*(sum(rr)+l)+ld] = var(str(y)+str(l)+str('_')+str(ld))
        for key in range(1,sum(rr)+1): 
            for l in range(1,d+1): z[key*d+l] = var(str(x)+str(key)+str(x)+str(l))
        edgesrhoG = [i for n, i in enumerate(rhoG) if i not in rhoG[:n]]
        vertrhoG = set(flatten(edgesrhoG));
        m=len(EP); 
        for l in range(m): vertrhoG.remove(sum(rr)+l+1);
        strr = '*λ'*len(vertrhoG)
        for i in vertrhoG:
            for l in range(1,d+1): strr = '*mu({},{},{})'.format(z[i*d+l], λ, β) + strr
            for l in range(1,d+1): strr = strr + ').integrate({},-infinity,+infinity)'.format(z[i*d+l])
        for i in edgesrhoG:
            for l in range(1,d+1): strr = '*H({},{},{})'.format(z[i[0]*d+l],z[i[1]*d+l], β) + strr
        strr = '('*len(vertrhoG)*d+strr[1:]
        jcumulants += eval(preparse(strr))
    print("\n");
    jcumulants = simplify(jcumulants).canonicalize_radical().maxima_methods().rootscontract().simplify()
    return jcumulants._sympy_()
\end{lstlisting} 

\subsubsection*{Acknowledgements}
\noindent
We thank Xueying Yang for essential contributions to the
SageMath cumulant codes. 

\footnotesize

\newcommand{\etalchar}[1]{$^{#1}$}
\def\cprime{$'$} \def\polhk#1{\setbox0=\hbox{#1}{\ooalign{\hidewidth
  \lower1.5ex\hbox{`}\hidewidth\crcr\unhbox0}}}
  \def\polhk#1{\setbox0=\hbox{#1}{\ooalign{\hidewidth
  \lower1.5ex\hbox{`}\hidewidth\crcr\unhbox0}}} \def\cprime{$'$}

\end{document}